\documentclass[11pt]{article}
\input epsf
\usepackage{enumerate, fullpage}
\usepackage{graphicx}
\usepackage{multirow}
\usepackage{color}

\def\b1{\mathbf{1}}

\def\bM{\mathbf{M}}

\newcommand{\beq}{\begin{equation}}
\newcommand{\eeq}{\end{equation}}
\newcommand{\bmath}[1]{{\boldmath ${#1}$}}

\newcommand{\beqnr}{\begin{eqnarray}}
\newcommand{\eeqnr}{\end{eqnarray}}

\newcommand{\benum}{\begin{enumerate}}
\newcommand{\eenum}{\end{enumerate}}

\newcommand{\argmin}{\mathop{\rm argmin}}
\newcommand{\QED}{\rule{.1in}{.1in}}
\newcommand{\cA}{{\cal A}}
\newcommand{\cB}{{\cal B}}
\newcommand{\cC}{{\cal C}}

\newcommand{\cE}{{\cal E}}
\newcommand{\cF}{{\cal F}}
\newcommand{\cG}{{\cal G}}

\newcommand{\cD}{{\cal D}}
\newcommand{\cR}{{\cal R}}

\newcommand{\cN}{{\cal N}}

\newtheorem{DE}{Definition}[section]

\newtheorem{LE}[DE]{Lemma}
\newtheorem{EX}[DE]{Example}
\newtheorem{RE}[DE]{Remark}
\newtheorem{THM}[DE]{Theorem}

\newcommand{\qed}{\mbox{}\hspace*{\fill}\nolinebreak\mbox{$\rule{0.7em}{0.7em}$}}

\advance \topmargin by -\headheight
\advance \topmargin by -\headsep

\evensidemargin \oddsidemargin

\begin{document}

\begin{center}

\Large {\bf The N - k Problem in Power Grids:  New Models, 
Formulations and Numerical Experiments \\ (extended version)}\footnote{Partially funded by award NSF-0521741 and award DE-SC000267}

\vskip 0.5cm

\Large

Daniel Bienstock
 and Abhinav Verma\\
Columbia University \\
New York 
\vskip 0.5cm
{\Large May 2008; Revised December 2009}
\end{center}

\begin{abstract}
Given a power grid modeled by a network together with equations describing
the power flows, power generation and consumption, and the laws of physics,
the so-called $N - k$ problem asks whether there exists a set of $k$ or fewer
arcs whose removal will cause the system to fail.  The
case where $k$ is small is of practical interest. We present 
theoretical and computational results involving a mixed-integer
model and a continuous nonlinear model related to this question.
\end{abstract}

\section{Introduction}\label{intro}

Recent large-scale power grid failures have highlighted the need for effective
computational tools for analyzing vulnerabilities of electrical 
transmission networks.  Blackouts are extremely rare, but their 
consequences can be severe.   Recent blackouts had, as their root
cause, an exogenous damaging event (such as a storm) which 
developed into a system collapse even though the initial quantity of 
disabled power lines  was small.  

As a result, a problem that has gathered 
increasing importance is what might be termed the {\em vulnerability evaluation} problem:
given a power grid, is there a small set of power lines whose removal will
lead to system failure?  Here, ``smallness'' is parametrized by an integer
$k$, and indeed experts have called for small values of $k$ (such as $k = 3$ or
$4$) in the analysis.  Additionally, an explicit goal in the formulation of
the problem is that the analysis should be agnostic: we are interested in
rooting out small, ``hidden'' vulnerabilities of a complex system which is
otherwise quite robust;  as much as possible the search for such vulnerabilities
should be devoid of assumptions regarding their structure.  

This problem is not new, and researchers have used a variety of names for it:
network {\em interdiction}, network {\em inhibition} and so on, although
the ``N - k problem'' terminology is common in the industry.  Strictly
speaking, one could consider the failure of different types of components (generators, transformers) in addition to lines, however, in this paper ``N'' 
will represent the number of lines.  We will 
provide a more complete review of the literature later
on;  the core central theme is that the $N - k$ problem is intractable, even
for small values of $k$ -- the pure enumeration approach is simply impractical.  
In addition to the combinatorial explosion, another
significant 
difficulty is the need to model the laws of physics governing power flows in a
sufficiently accurate and yet computationally tractable manner:  power flows
are much more complex than ``flows'' in traditional computer science or
operations research applications.

A critique that has been leveled against optimization-based approaches to the
$N - k$ problem is that they tend to focus on large values of $k$, say $k = 8$.  In our experience, when
$k$ is large the problem tends to become easier (possibly because many solutions exist) but on the other hand the argument
can be made that the cardinality of the attack is unrealistically large.  At the other
end of the spectrum lies the case $k = 1$, which can be addressed by enumeration but
may not yield useful information.  The middle range, $2 \le k \le 5$, is both relevant
and difficult, and is our primary focus.

We present
results using two optimization models.  The first (Section \ref{mip}) 
is a new linear mixed-integer programming
formulation that explicitly models a ``game'' between a fictional attacker
seeking to disable the network, and a controller who tries to prevent a 
collapse by selecting which generators to operate and
adjusting generator outputs and demand levels.  As far as we can tell, the 
problem we consider here is more general than has been previously studied in
the literature; nevertheless our approach yields practicable solution times
for larger instances than previously studied.

The second
model (Section \ref{nlp}) is given by a new, continuous 
nonlinear programming formulation whose goal is
to capture, in a compact way, the interaction between the underlying physics
and the network structure.  While both formulations provide substantial
savings over the pure enumerational approach, 
the second formulation appears particularly effective and scalable; enabling us
to handle in an optimization framework 
models an order of magnitude larger than those we have seen in the literature.

\subsubsection{Previous work on vulnerability problems}\label{prevwork}

There is a large amount of prior work on optimization methods applied to blackout-related problems.  \cite{pinar} includes a fairly comprehensive
survey of recent work.  

Typically work has focused on identifying a small set of arcs whose 
removal (to model complete failure) will result in a network unable to deliver
a minimum amount of demand.  A problem of this type can be 
solved using mixed-integer programming techniques
techniques, see \cite{alvarez}, \cite{swb}, 
\cite{arrgal05}.  We will review this work in
more detail below (Section \ref{prevwork2}).  Generally speaking, the mixed-integer programs to be solved
can prove quite challenging.

A different line of research on vulnerability problems focuses on attacks with certain structural properties, see \cite{bier}, \cite{pinar}.  
An example of this
approach is used in \cite{pinar}, where as an approximation to the 
$N - k$ problem with AC power flows, 
a linear mixed-integer program to solve the following combinatorial problem: remove a
minimum number of arcs, such that in the resulting network there is a
partition of the nodes into two sets, $N_1$ and $N_2$, such that
\begin{eqnarray}
&& D(N_1) \, + \, G(N_2) \, + \, cap( N_1, N_2 ) \, \le \, Q^{min}.\label{haha}
\end{eqnarray}
\noindent Here $D(N_1)$ is the total demand in $N_1$, $G(N_2)$ is the
total generation capacity in $N_2$, $cap( N_1, N_2 )$ is the total capacity
in the (non-removed) arcs between $N_1$ and $N_2$, and $Q^{min}$ is a minimum  
amount of demand that needs to be satisfied.  The quantity in the left-hand side in the
above expression is an upper-bound on the total amount of demand that
can be satisfied -- the upper-bound can be strict because under
power flow laws it may not be attained.

Thus this is an approximate model that could underestimate the effect of
an attack (i.e. the algorithm may produce attacks larger than strictly
necessary).  On the other hand, 
methods of this type bring to bear powerful 
mathematical tools, and
thus can handle larger problems than algorithms that rely on generic mixed-integer programming techniques.   
Our method in Section \ref{nlp} can also be
viewed as an example of this approach.\\

Finally, we mention that the most sophisticated models for the behavior of a 
grid under stress
attempt to capture the multistage nature of blackouts, and are thus
more comprehensive than the static models considered above and in this paper.
See, for example, \cite{CarrerasCH02}-\cite{CarrerasCAS04}, and
\cite{bie}.\

\subsubsection{Power Flows} \label{poflo}

Here we provide a brief introduction to the 
so-called {\em linearized}, or {\em DC} power flow model.  For the purposes of our
problem, a grid is represented by a directed network $\cN$, 
 where:
\begin{itemize}
\item Each node corresponds to a ``generator'' (i.e., a supply node), or to a ``load'' (i.e., a demand node), or to a node that neither generates nor 
consumes power.  We denote by $\cG$ the set of generator nodes.
\item If node $i$ corresponds to a generator, then there are values 
$0 \le P_i^{min} \le P_i^{max}$.   If the generator is operated, then 
its output must be in the range $[P_i^{min} , P_i^{max}]$; if the generator
is not operated, then its output is zero.  In general, we expect $P^{min}_i > 0$.

\item If node $i$ corresponds to a demand, then there is a value $D_i^{nom}$
(the ``nominal'' demand value at node $i$).
We will denote the set of demands by $\cD$.
\item The arcs of $\cN$ represent power lines.  For each arc $(i,j)$, we are given a parameter $x_{ij} > 0$ (the
resistance) and a parameter $u_{ij}$ (the capacity).
\end{itemize}

\noindent Given a set $\cC$ of operating generators, a {\em power flow} is a solution to the system
of constraints given next.  In this system, for each
   arc $(i, j)$, we use a variable $f_{ij}$ to represent the (power) flow on 
$(i, j)$ -- possibly $f_{ij} < 0$, in which case power is 
effectively flowing from $j$ to $i$.  In addition, for each node $i$ we
will have a variable $\theta_i$ (the ``phase angle'' at $i$).  Finally, if $i$ is a generator node,
then we will have a variable $P_i$, while if $i$ represents a demand
node, we will have a variable  $D_i$.  Given a node $i$, we represent with 
$\delta^{+}(i)$ ($\delta^{-}(i)$) the set of arcs oriented out of (respectively, into) $i$.\\

\noindent In the system given below, hereafter denoted by \bmath{P(\cN, \cC)},
constraints (\ref{1b}), (\ref{1f}) and (\ref{1g}) are typical for network flow models  (for
background see \cite{amo}) and model, respectively: 
flow balance (i.e., net flow leaving a node equals net supply at that node),
generator and demand node bounds. Constraint (\ref{ohm-eq}) is a commonly used linearization of 
more complex equations describing power flow physics; for background see \cite{lecture}.   We will comment on this equation in Section \ref{acdc}.

\begin{eqnarray}
& & \sum_{(i,j)\in \delta^{+}(i)} f_{ij} - \sum_{(j,i)\in \delta^{-}(i)} f_{ji}  = \left \{
\begin{array}{lll}
~~P_i & {i \in \cC} &\\
-D_i & {i \in \cD} &  \\
~~0   & \mbox{otherwise} &
\end{array} \right. \label{1b}\\
&& \nonumber \\
 & &  \theta_i - \theta_j - x_{ij}f_{ij} = 0  \ \ \ \forall (i,j) \label{ohm-eq}\\
&& \nonumber \\
 & &  | f_{ij} |  \ \ \ \le \ \ \ u_{ij}, \ \ \ \forall (i,j) \label{1ap}\\
&& \nonumber \\
&&   P_i^{min} \leq P_i \leq P_i^{max} \ \ \ \forall i \in \cC \label{1f}\\
&& \nonumber \\
&& 0 \leq D_j \leq D_j^{nom} \ \ \ \forall j \in \cD \label{1g}
\end{eqnarray}
\noindent One could also impose explicit upper bounds on the quantities
$|\theta_i - \theta_j|$ (over the arcs $(i,j)$).  However, note that our model already
does so, implicitly: because of constraint (\ref{ohm-eq}), 
imposing an upper bound on $| f_{ij}|$ (i.e., eq. (\ref{1ap}))  
is equivalent to imposing an upper bound on $|\theta_i - \theta_j|$.

\subsubsection{Basic properties}
A useful property satisfied by the linearized model is given by the 
following result.

\begin{LE} \label{unique}
Let $\cC$ be given, and suppose $\cN$ is connected.  Then 
for any choice of nonnegative values $P_i$ (for $i \in \cC$) and
$D_i$ (for $i \in \cD$) such that 
\begin{eqnarray}
\sum_{i \in \cC} P_i & = & \sum_{i \in \cD} D_i, \label{suppdembal}
\end{eqnarray}
\noindent system (\ref{1b})-(\ref{ohm-eq}) has a unique solution in the $f_{ij}$; 
thus, 
the solution is also unique in the
$\theta_i - \theta_j$ (over the arcs $(i,j)$).
\end{LE}
\noindent {\em Proof.} Let $N$ denote the node-arc incidence of the network \cite{amo}, let
$b$ be the vector with an entry for each node, where $b_i = P_i$ for $i \in \cC$, $b_i = -D_i$ 
for $i \in \cD$, and $b_i = 0$ otherwise.  Writing $X$ for the diagonal matrix with entries
$x_{ij}$, then (\ref{1b})-(\ref{ohm-eq}) can be summarized as
\begin{eqnarray}
N f & = & b, \label{fc}\\
N^T \theta \, - \, X f & = & 0. \label{fo}
\end{eqnarray}
\noindent Pick an arbitrary node $v$; then system (\ref{fc})-(\ref{fo}) has a solution iff it has 
one with $\theta_v = 0$.  As is well-known, $N$ does not have full row rank, 
but writing $\bar N$ for the submatrix of $N$ with the row corresponding
to $v$ omitted then the connectivity assumption implies that $\bar N$ does have full row rank
\cite{amo}.
In summary, writing $\bar b$ for the corresponding subvector of $b$, we have that 
(\ref{fc})-(\ref{fo}) has a solution iff 
\begin{eqnarray}
\bar N f & = & \bar b, \label{fc2} \\
\bar N^T \eta \, - \, X f & = & 0. \label{fo2}
\end{eqnarray}
\noindent where the vector $\eta$ has an entry for every node other than $v$.  Here,
(\ref{fo2}) implies $f = X^{-1} \bar N^T \eta$, and so (\ref{fc2}) implies that 
$\eta = (\bar N X^{-1} \bar N^T)^{-1} \bar b$ (where the matrix is invertible since $\bar N$ has
full row rank). Consequently, $f = X^{-1} \bar N^T(\bar N X^{-1} \bar N^T)^{-1} \bar b$. \QED

\begin{RE} We stress that Lemma \ref{unique} concerns the subsystem of \bmath{P(\cN, \cC)} consisting of 
(\ref{1b}) and (\ref{ohm-eq}). In particular, the ``capacities'' $u_{ij}$ 
play no role in the determination of solutions.  
\end{RE}
\noindent When the network is not connected Lemma \ref{unique} can be extended
by requiring that (\ref{suppdembal}) hold for each component.

\begin{DE}
Let $(f, \theta, P, D)$ be feasible a solution to $P(\cN, \cC)$. The {\bf throughput} of $(f, \theta, P, D)$ is defined as
\begin{eqnarray}
\frac{ \sum_{i \in \cD} D_i }{ \sum_{i \in \cD} D^{nom}_i }.
\end{eqnarray}
\noindent The throughput of $\cN$ is the maximum throughput of any feasible
solution to $P(\cN, \cC)$.
\end{DE}

\subsubsection{DC and AC power flows, and other modeling issues}\label{acdc}

Constraint (\ref{ohm-eq}) is reminiscent of Ohm's equation --
in a direct current (DC) network (\ref{ohm-eq}) precisely represents Ohm's equation. In the case of an AC network (the relevant case when dealing with power
grids) (\ref{ohm-eq}) only {\em approximates} a complex system
of nonlinear equations (see \cite{lecture}).  The issue of whether to use
the more accurate nonlinear formulation, or the approximate DC formulation,
is rather thorny.  On the one hand, the linearized
formulation certainly is an approximation only.  On the other hand,
a formulation that models AC power flows can prove intractable or
may reflect difficulties inherent with the underlying real-life problem
(e.g. the formulation may have multiple solutions).  

We can provide a very brief summary of how the linearized model arises.  
First, AC power flow models typically include equations of the form
\begin{eqnarray}
x_{ij} f_{ij} & = & \sin(\theta_i - \theta_j), \ \ \ \forall (i,j), \label{sineq}
\end{eqnarray}
\noindent as opposed to (\ref{ohm-eq}). Here, the $f$ quantities describe ``active'' power flows and the $\theta$ describe 
phase angles.   In normal operation of a transmission system, one would expect that
$\theta_i \approx \theta_j$ for any arc $(i,j)$ and thus (\ref{sineq}) can be linearized
to yield (\ref{ohm-eq}).  Hence the linearization is only valid if we additionally impose that
$| \theta_i - \theta_j |$ be {\em very} small.  However, in the literature
one sometimes sees this ``very small'' constraint relaxed, for the reason that we are interested
in regimes where the network is {\em not} in a normal operative mode.  One might still 
impose {\em some} explicit upper bound on $| \theta_i - \theta_j |$; we have seen publications with and
without such a constraint.  As explained above, our model already includes
such a constraint in implicit form.  

In either case, one ends up using a representation of the 
problem that is arguably {\em not} valid. But the key 
fact is that constraint (\ref{sineq}) gives rise to extremely 
complex models.  Studies that require multiple power flow computations tend
to rely on the linearized formulation, in the expectation that some
useful information will emerge.  This will be the approach we
take in this paper, though some of our techniques extend directly to 
AC models and this will remain a topic for future research.

An approach such
as ours can therefore be criticized because it relies on an ostensibly 
approximate model; on the other hand we are able to focus more explicitly on
the basic combinatorial complexity that underlies the $N-k$ problem.  In
contrast, an approach addressing AC model characteristics would have a better representation of the physics, but at the cost of not being able to tackle the combinatorial
complexity quite as effectively, for the simple reason that the theory 
and computational machinery for linear programming are far more mature,
effective {\em and} scalable than those
for nonlinear, nonconvex optimization.  
In summary, both approaches present 
limitations and benefits.  In this paper, our bias is toward explicitly handling
the combinatorial nature of the problem.\\

\noindent A final point that we would like to stress is that whether we use an
AC or DC power flow model, the resulting problems  have a far more complex structure
than (say) traditional single- or multi-commodity flow models because of 
side-constraints  such as (\ref{ohm-eq}).  Constraints of this type give
rise to counter-intuitive behavior reminiscent of 
Braess's Paradox \cite{braess}.\\

\noindent The model we consider in the paper can be further enriched by including many 
other real-life constraints.  For example, when considering the response to a contingency one
could insist that demand be curtailed in a (geographically) even-handed pattern, and not 
in an aggregate fashion, as we do below and has been done in the literature.  Or we could 
impose upper bounds on the number of stand-by power units that are turned on in the event
of a contingency (and this itself could have a geographical perspective).  Many such realistic
features suggest themselves and could give rise to interesting extensions to the problem we
consider here.

\section{The ``N - k'' problem}\label{nmk}

Let $\cN$ be a network with $n$ nodes and $m$ arcs representing a power grid.  We denote
the set of arcs by $E$ and the set of nodes by $V$.
A fictional {\em attacker} 
wants to remove a small number of arcs
from $\cN$ so that in the resulting  network all feasible flows should have
low throughput.   At the same time, a {\em controller} is operating the network;  the controller responds to an attack by appropriately choosing the set
$\cC$ of operating generators, their output levels, and the demands $D_i$, so as to 
feasibly obtain high throughput.  Note that the 
controller's adjustment of demands
corresponds to the traditional notion of ``load shedding'' in the
power systems literature.

Thus, the attacker seeks to defeat {\em all possible courses of action} by
the controller, in other words, we are modeling the problem as a Stackelberg
game between the attacker and the controller, where the attacker moves first. To cast this in a precise way we 
will use the following definition. We 
let $0 \le T^{min} \le 1$ be a given value.

\begin{DE} \label{feas}
Given a network $\cN$,
\begin{itemize}
\item An {\bf attack} $\cA$ is a set of arcs removed by the attacker.  
\item Given an attack $\cA$, the {\bf surviving network} $\cN - \cA$ is
the subnetwork of $\cN$ consisting of the arcs not removed by the attacker.
\item A {\bf configuration} is a set $\cC$ of generators.
\item We say that an attack $\cA$ {\bf defeats} a configuration $\cC$, if either (a) 
the maximum throughput of any feasible solution to \bmath{P(\cN - \cA , \cC)}
is strictly less than 
$T^{min}$, or (b) no feasible solution to \bmath{P(\cN - \cA , \cC)} exists. 
Otherwise we say that $\cC$ defeats $\cA$.  

\item We say that an attack is {\bf successful}, if it defeats {\bf every}
configuration.
\item The {\bf min-cardinality attack problem} consists in finding a successful attack $\cA$ with $|\cA|$ minimum. 
\end{itemize}
\end{DE}

\noindent Our primary focus will be on the min-cardinality attack problem.  Before proceeding further we would like to comment on our model, 
specifically on the parameter $T^{min}$.  In a practical use
of the model, one would wish to experiment with different values
for $T^{min}$ -- for the simple reason that an attack $\cA$ which is not
successful for a given choice for $T^{min}$ could well be successful for a
slightly larger value; e.g. no attack of cardinality $3$ or less exists that 
reduces demand by $31 \%$, and yet there exists an attack of cardinality $3$ 
that reduces
satisfied demand by $30 \%$. In other words, the minimum cardinality of a successful attack could vary 
substantially as a function of $T^{min}$.

Given this fact, {\em it might appear} that a better approach to the
power grid vulnerability problem would be to leave out the parameter 
$T^{min}$ entirely, and instead
redefine the problem to that of finding a set of $k$ or fewer arcs 
to remove, so that the resulting network has minimum throughput (here,
$k$ is given).
We will refer to this as the 
{\em budget-k min-throughput problem}.   However, there are reasons
why this latter problem is less attractive than the min-cardinality problem.

\begin{itemize}
\item [(a)] The min-cardinality and
min-throughput problems are duals of each other. A modeler considering
the min-throughput problem would want to run that model multiple times, 
because given $k$, the value of the budget-$k$ min-throughput problem could be 
much smaller than the value of the budget-$(k+1)$ min-throughput problem.
For
example, it could be the case that using a budget of $k = 2$, the attacker
can reduce throughput by no more than $5 \%$; 
but nevertheless with a budget of $k = 3$,  
throughput can be reduced by e.g. $50 \%$.  
In other words, even if a network is ``resilient'' against attacks of size $\le 2$, it
might nevertheless prove very vulnerable to attacks of size $3$.  
For this reason, and given that the models
of grids, power flows, etc., are rather approximate, a practitioner would
want to test various values of $k$ -- this issue
is obviously related to what percentage of demand loss would be considered
tolerable, in other words, the parameter $T^{min}$.
\item[(b)] From an operational perspective it should be straightforward to
identify reasonable values for the quantity $T^{min}$; whereas the value 
$k$ is more obscure and bound to models of how much power the adversary can
wield.
\item [(c)] Because of a subtlety that arises from having positive quantities
$P^{min}_i$, explained next, it turns out that the min-throughput problem
is significantly more complex and is difficult to even formulate in a 
compact manner.
\end{itemize}

\noindent We will now expand on (c). One would
expect that when a configuration $\cC$ is defeated by an attack $\cA$, it is because only
small throughput solutions are feasible in $\cN-\cA$.
However, because the lower bounds
$P_i^{min}$ are in general strictly positive, it may also be the case that
{\em no feasible solution to} \bmath{P(\cN - \cA , \cC)} {\em exists}.    

\begin{EX}\label{subtle}

Consider the following example on a network $\cN$ with three nodes (see Figure \ref{fig:22}), where
\begin{enumerate}
\item Nodes $1$ and $2$ represent generators; $P^{min}_1 = 2$, $P^{max}_1 = 4$, $P^{min}_2 = 0$, and $P^{max}_2 = 4$, 
\item Node $3$ is a demand node with $D^{nom}_3 = 6$.  Furthermore, $T^{min} = 1/2$.
\item There are three arcs; arc $(1,2)$ with $x_{12} = 1$ and $u_{12} = 1$, arc $(2,3)$ with $x_{23} = 1$ and $u_{23} = 5$, 
and arc $(1,3)$ with $x_{13} = 1$ and $u_{13} = 3$. 
\end{enumerate}
\begin{figure}[h] 
\centering
\includegraphics[width=3in,angle=270]{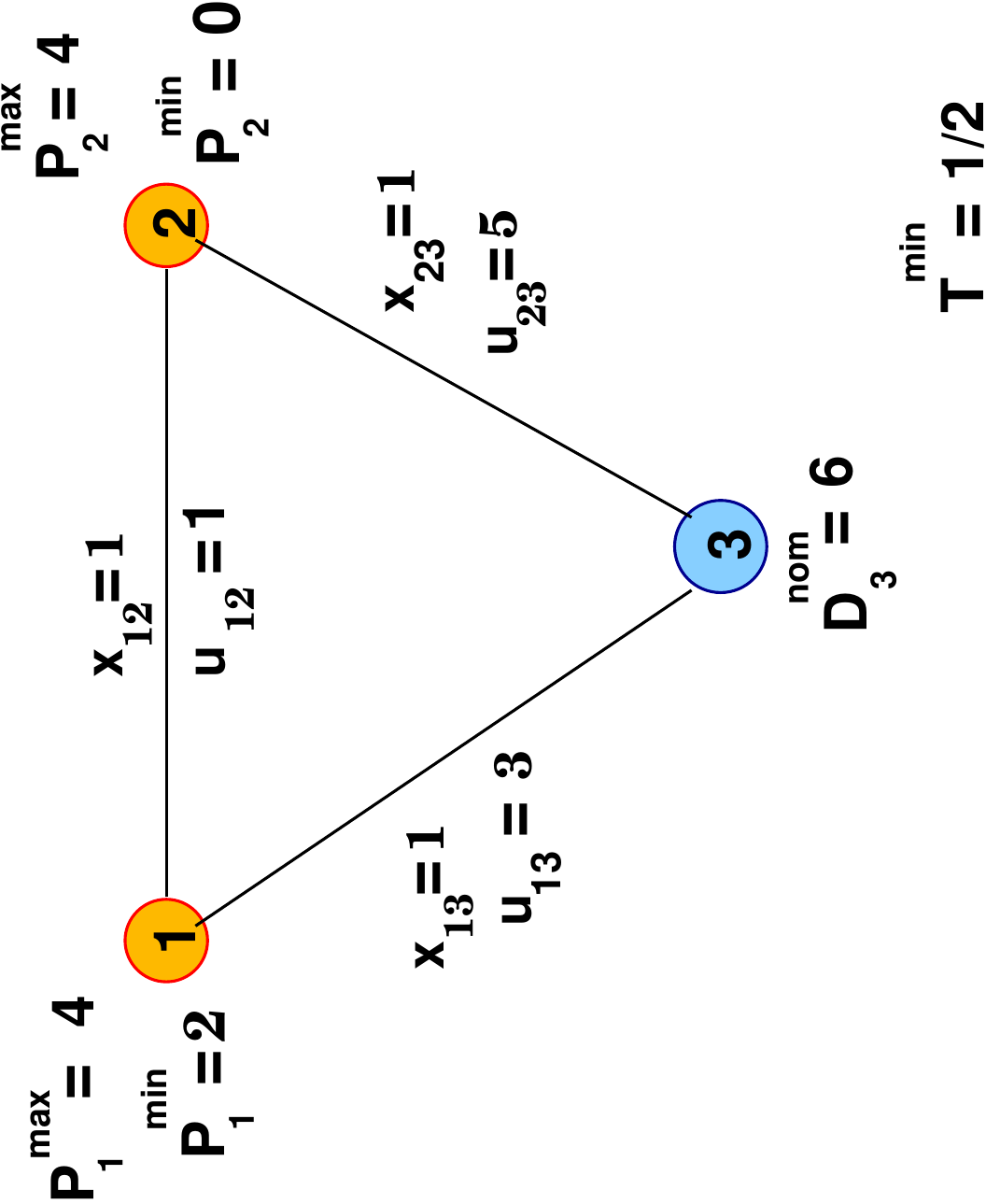}
\caption{A simple example where the controller must turn off a generator.\label{fig:22}}
\end{figure}

\noindent When the network is not attacked, the following solution is feasible: $P_1 = P_2 = 3$, $D_3 = 6$, $f_{12} = 0$, 
$f_{13} = f_{23} = 3$, $\theta_1 = \theta_2 = 0$, $\theta_3 = -3$.  This solution has throughput $100\%$.  On the other hand,
consider the attack $\cA$ consisting of the single arc $(1,3)$. 
Since $P^{min}_1 > u_{12}$, \bmath{P(\cN - \cA , \cC)} has no feasible solution, and $\cA$ defeats $\cC$, for $\cC = \{ 1 , 2 \}$ or $\cC = \{1\}$,  in spite of the fact that that total generator capacity is
enough to meet at least $2/3 > T^{min}$ of the demand.
Thus,
$\cA$ is successful if and only if it also defeats the configuration where we only operate generator $2$, which 
it does not since in that configuration we can feasibly send up to four units of flow on $(2,3)$.  \end{EX}

As the example highlights, it is important to understand how 
an attack $\cA$ can defeat a particular configuration $\cC$. It
turns out that there are {\em three} different ways for this to happen.  
\begin{itemize}
\item[(i)] Consider a partition of the nodes of $\cN$ into two classes, 
$N^1$ and $N^2$.  Write
\begin{eqnarray}
D^k & = & \sum_{i \in \cD \cap N^k} D^{nom}_i, \,\,\,\,k = 1, 2, \,\,\,\,\,\, \mbox{and} \\
P^k & = & \sum_{i \in \cC \cap N^k} P^{max}_i, \,\,\,\,k = 1, 2,
\end{eqnarray}
\noindent e.g. the total (nominal) demand in $N_1$ and $N_2$, and the 
 maximum power generation in $N_1$ and $N_2$, respectively.  The following
condition, should it hold, would guarantee that $\cA$ defeats $\cC$:
\begin{eqnarray}
T^{min} \sum_{j \in \cD} D^{nom}_j - \min\{ D^1, P^1\} - \min\{D^2, P^2\} & > & \sum_{(i,j) \notin \cA \, : \, i \in N^1, \, j \in N^2} u_{ij} \, + \,  \nonumber \\
&& \,\,\,\,\, \sum_{(i,j) \notin \cA \, : \, i \in N^2, \, j \in N^j} u_{ij}. \label{cond1}
\end{eqnarray}
\noindent To understand this condition, note that for $k = 1, 2$, 
$\min\{ D^k, P^k\}$ is the maximum demand within $N^k$ that could possibly 
be met using power flows that do not leave $N^k$. Consequently 
the left-hand side of (\ref{cond1}) is a lower bound on the amount of flow
that must travel between $N^1$ and $N^2$, whereas the right-hand side 
of (\ref{cond1}) is the total capacity of arcs between $N^1$ and $N^2$ 
under attack $\cA$.  In other words, condition (\ref{cond1}) amounts to a 
mismatch between demand and supply.  A special case of (\ref{cond1}) is that
where in $\cN - \cA$ there are no arcs between $N^1$ and $N^2$, i.e. the
right-hand side of (\ref{cond1}) is zero. Condition (\ref{cond2}) is
similar, but not identical, to (\ref{haha}).
\item [(ii)] Consider a partition of the nodes of $\cN$ into two classes, 
$N^1$ and $N^2$.  Then attack $\cA$ defeats $\cC$ if
\begin{eqnarray}
\sum_{i \in \cD \cap \in N^1} D^{nom}_i \ + \ \sum_{(i,j) \notin \cA \, : \, i \in N^1, \, j \in N^2} u_{ij} & < & \sum_{i \in \cC \cap N^1} P^{min}_i, \label{cond2}
\end{eqnarray}
\noindent i.e., the minimum power output within $N^1$ exceeds the
maximum demand within $N^1$ plus the sum of arc capacities leaving $N^1$.  Note that (ii) may apply even if (i) does not.  
\item [(iii)] Even if (i) and (ii) do not hold, it may still be the case 
that the system (\ref{1b})-(\ref{1g}) does not admit a feasible solution.
To put it differently, suppose that for every choice of demand values $0 \le D_i \le D^{nom}_i$ (for $i \in \cD$)  and supply values $P^{min}_ i \le P_i \le P^{max}_i$ 
(for $i \in \cC$) such that 
$\sum_{i \in \cC} P_i \, = \, \sum_{i \in \cD} D_i$ the
unique solution to system (\ref{1b})-(\ref{ohm-eq}) in network $\cN - \cA$
(as per Lemma \ref{unique}) does {\em not} satisfy the ``capacity'' 
inequalities $| f_{ij} |  \le u_{ij}$ for all arcs $(i,j) \in \cN - \cA$. Then
$\cA$ defeats $\cC$.  This is the most subtle case of all -- it involves
the interplay of constraints (\ref{ohm-eq}) and (\ref{1b}) in the power flow
model.
\end{itemize}

Note that in particular in case (ii), the defeat condition is unrelated to
throughput.  Nevertheless, should case (ii) arise, it is clear that
the attack has succeeded (against configuration $\cC$) -- this 
makes the min-throughput problem difficult to model;
our formulation for the min-cardinality
problem, given in  Section \ref{mip}, does capture the
three defeat criteria above.  \\

From a practical perspective, it is important to handle models where
the values $P^{min}_i$ are positive. It is also
important to model {\em standby} generators that are turned on when needed, 
and to model the turning off of generators that are unable to dispose of their
minimum power output, post-attack.  All these features arise in 
practice.  
Example \ref{subtle} above shows that models where generators cannot be turned off
can exhibit unreasonable behavior.  Of course, the ability to select the 
operating generators comes
at a cost, in that in order to certify that an attack is successful we need
to evaluate, at least implicitly, a possibly exponential number of control possibilities.\\

As far as we can tell, most (or all) prior work in the literature {\bf does} require that
the controller must always use the configuration $\bar \cG$ consisting of all generators.  As the example shows, however,
if the quantities $P^{min}_i$ are positive there may be attacks $\cA$ such 
that \bmath{P(\cN - \cA , \bar \cG)} is infeasible.  Because of this fact, algorithms that rely on direct application of
Benders' decomposition or bilevel programming are problematic, and
{\em invalid} formulations can be found in the literature.   \\

Our approach works with general $P^{min} \ge 0$ quantities; thus, it also applies to the case
where we always have $P^{min}_i = 0$. In this case our formulation is simple enough that 
a commercial integer program solver can directly handle instances larger than previously
reported in the literature.

\subsubsection{Non-monotonicity in optimal attacks}

\begin{figure}[htb] 
\centering
\includegraphics[width=6in]{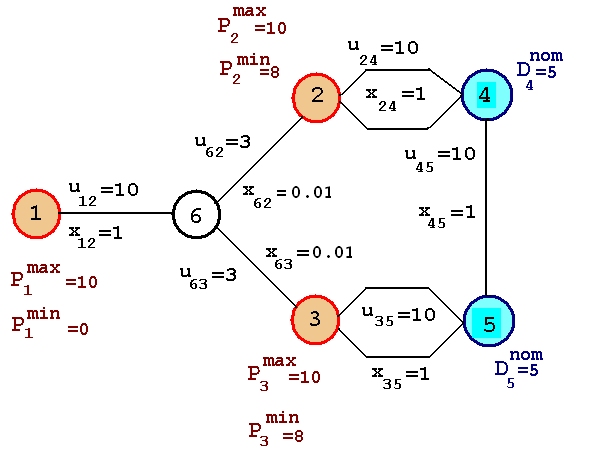}
\caption{An example where there is a successful attack with $k = 1$ but none with $k = 2$.\label{fig:nonmon}}
\end{figure}

Consider the example in Figure \ref{fig:nonmon}, where we assume $T^{min} = 0.3$.  Notice that there are two parallel copies of arcs $(2,4)$ and $(3,5)$, each with capacity $10$ and
impedance $1$.  
It is
easy to see that the network with no attack is feasible: we operate generator $1$ and not
operate generators $2$ and $3$, and send $3$ units of flow along the paths 
$1-6-2-4$ and $1-6-3-5$ (the flow on e.g. the two parallel $(2,4)$ arcs is
evenly split).

On the other hand, consider the attack consisting of arc $(1,6)$ -- 
we will show this attack
is successful.  To see this, note that under this attack, 
the controller cannot operate both generators
$2$ and $3$, since their combined minimum output exceeds the total demand.  Thus, without loss of generality suppose
that only generator $3$ is operated, and assume by contradiction that a feasible solution
exists -- then this solution must route {\em at most} $3$ units of flow along
$3-6-2-4$, and (since $P^{min}_3 = 8$) {\em at least} $5$ units of flow on $(3,5)$ 
(both copies altogether).  In such
a case, the phase angle drop from $3$ to $5$ is at least $2.5$, whereas the drop from
$3$ to $4$ is at most $1.56$.  In other words, $\theta_4 - \theta_5 \ge 0.94$, and so
we will have $f_{45} \ge 0.94$  -- thus, the net inflow at node 
$5$ is at least $5.94$.  Hence the attack is indeed successful.  

However, there is {\em no} successful attack consisting 
of arc $(1,6)$ and another arc.  To see this, note that if one of $(2,6)$, 
$(3,6)$ or $(4,5)$ are also removed then the controller can {\em just} 
operate one of the two generators $2$, $3$ and meet eight units of demand. 
Suppose that (say) one of the two copies of $(3,5)$ is removed (again, in 
addition to $(1,6)$).  Then the controller operates generator $2$, sending
$2.5$ units of flow on each of the two parallel $(2,4)$ arcs; 
thus $\theta_2 - \theta_4 = 2.5$. The controller also routes $3$ units of flow 
along $2-6-3-5$, and therefore $\theta_2 - \theta_5 = 3.06$. Consequently
$\theta_4 - \theta_5 = .56$, and $f_{45} = .56$, resulting in a feasible flow which 
satisfies $4.44$ units of demand at $4$ and $3.56$ units of demand at $5$. 

In fact, it can be shown that {\em no} successful attack 
of cardinality $2$ exists -- hence we observe non-monotonicity.

By elaborating on the above, one can create 
examples with arbitrary patterns in the cardinality of successful attacks.
One can also generate examples that exhibit non-monotone behavior in response to
controller actions.  In both cases, the non-monotonicity can be viewed as a manifestation
of the so-called ``Braess's Paradox'' \cite{braess}.  In the above example 
we can observe combinatorial subtleties that arise from the ability of the 
controller to choose which generators to operate, 
and from the lower bounds on output in operating generators. Nevertheless,
it is clear that the critical core reason for the 
complexity is the interaction between 
phase angles and flows, i.e. between ``Ohm's law'' (\ref{ohm-eq}) and 
flow conservation (\ref{1b}) -- the combinatorial attributes of the problem 
exercise this interaction. 

\subsubsection{Brief review of previous work on related problems}\label{prevwork2}

The min-cardinality problem, as defined above,
can be viewed as a bilevel 
program where both the master problem and the subproblem 
are mixed-integer programs -- the master problem corresponds
to the attacker (who chooses the arcs to remove) and the subproblem to the
controller (who chooses the generators to operate).  
In general, such problems are
extremely challenging.  A recent general-purpose algorithm for such
integer programs is given in \cite{dnralphs}.

Alternatively, each configuration of generators can be viewed as a 
``scenario''.  In this sense our problem resembles a stochastic program, 
although without a probability distribution.  Recent work \cite{lind} 
considers a single commodity max-flow problem 
under attack by an interdictor with a limited attack budget; 
where an attacked arc 
is removed probabilistically, leading to a stochastic program (to minimize
the expected max flow). 
A deterministic, multi-commodity version of the same problem is given in
\cite{limsmith}.  \\

Previous work on the power grid vulnerability models proper has 
focused on cases where either 
the generator lower bounds $P^{min}_i$ are all zero,
or all generators must be operated (the single configuration case).  Algorithms for these problems have either relied on
heuristics, or on mixed-integer programming techniques, usually a direct use of Benders' decomposition or bilevel programming.  
\cite{alvarez} considers a version of the min-throughput problem with
$P^{min}_i = 0$ for all generators $i$, and presents an algorithm using
Benders' decomposition (also see references therein).  They analyze the so-called
IEEE One-Area and IEEE Two-Area test cases, with, respectively, 24 nodes and 38 arcs, and 48 
nodes and 79 arcs.  Also see \cite{swb}.

\cite{arrgal05} studies the IEEE One-Area test case, and allows $P^{min}_i > 0$, but does not
allow generators to be turned off; the authors present a bilevel programming
formulation which unfortunately is incorrect, due to reasons outlined above.\\ 

\subsection{A mixed-integer programming algorithm for the min-cardinality problem}\label{mip}
In this section we will describe an iterative algorithm for the min-cardinality
attack problem.  The algorithm iterates in Benders-like fashion, solving
at each iteration two mixed-integer programs.  Before describing the algorithm
we need to introduce some notation and concepts.

Let $\cA$ be a given attack.  Suppose the controller attempts to defeat
the attacker by choosing a certain configuration $\cC$ of generators.
Denote by $z^\cA$ the indicator vector for $\cA$, i.e. $z^\cA_{ij} = 1$ iff
$(i,j) \in \cA$. The controller needs to operate the network so as to
satisfy the required amount of demand without exceeding any arc capacity.
The last requirement can be cast as an optimization goal: minimize the
maximum arc overload subject to satisfying the desired level of demand.
Should this min-max overload be strictly greater than 1, the controller
is defeated.  In summary, the controller needs to solve the following linear program (where the variable ``t'' indicates the min-max overload):
\begin{eqnarray}
\mbox{\bmath{K_{\cC} (\cA):}} \hspace{.6in} t_{\cC}(z^\cA) && \doteq \ \mbox{ \ min} \ t \label{mintobj}\\
&& \mbox{Subject to:}  \nonumber \\
& & \sum_{(i,j)\in \delta^{+}(i)} f_{ij} - \sum_{(j,i)\in \delta^{-}(i)} f_{ji}  = \left \{
\begin{array}{lll}
~~P_i & {i \in \cG} &\\
-D_i & {i \in \cD} &  \\
~~~0   & \mbox{otherwise} &
\end{array} \right. \label{2b}\\
&& \nonumber \\
& &  \theta_i - \theta_j - x_{ij}f_{ij} = 0  \ \ \ \forall \ (i,j)  \notin \cA \label{2a}\\
&& \nonumber \\
 & &  u_{ij}~ t \ - \ |f_{ij}|  \ \ \ \ge \ \ \ 0, \ \ \ \forall \ (i,j) \notin \cA \label{2ap}\\
&& \nonumber \\
 & &  f_{ij}  \ \ = \ \ \ 0, \ \ \forall (i,j) \in \cA \label{2aq}\\
&& \nonumber \\
&&   P_i^{min} \ \leq P_i \ \leq P_i^{max} \ \ \ \forall i \in \cC \label{2f}\\
&& \nonumber \\
&& P_i \ \ = \ \ 0, \ \ \ \forall i \in \cG \ - \ \cC \label{2fp}\\
&& \nonumber \\
&& \sum_{j \in \cD} D_j ~ \geq ~ T^{min} \left( \sum_{j \in \cD} D^{nom}_j \right), \label{2pg} \\
&& \nonumber \\
&& 0 \leq ~ D_j ~ \leq D_j^{nom} \ \ \ \forall j \in \cD \label{2g}
\end{eqnarray}
\begin{RE} \label{infdef} Using the convention that the value of an infeasible linear 
program is infinite, $\cA$ defeats $\cC$ if and only if
$t_{\cC}(z^\cA) > 1$.
\end{RE}

Thus, an attack $\cA$ is {\em not} successful if and only if we can find
$\cC \subseteq \cG$ with 
$t_{\cC}(z^\cA) \le 1$; we test for this condition
by solving the problem:
$$ \min_{\cC \subseteq \cG} ~ t_{\cC}(z^\cA).$$
\noindent  This is done by replacing, in the above formulation, equations (\ref{2f}), (\ref{2fp}) with 
\begin{eqnarray}
&&  P_i^{min} y_i \ \leq P_i \ \leq P_i^{max} y_i,\ \ \ \forall i \in \cG, \label{2h} \\
&& y_i \ \ = 0 \ \ \mbox{or} \ \ 1, \ \ \ \forall i \in \cG. \label{2i}
\end{eqnarray}
\noindent Here, $y_i = 1$ if the controller operates 
generator $i$. \\

\noindent The min-cardinality attack problem can now be written as follows:
\begin{eqnarray}
&&  \mbox{ \ min} \ \sum_{(i,j)} z_{ij} \label{mincard} \\
&& \nonumber \\
&& t_{\cC}(z) \, > \, 1, \,\,\,\,\,\, \forall \,\, \cC \subseteq \cG, \label{nonlinear} \\
&& \nonumber \\
&& z_{ij} \, = \, 0 \,\, \mbox{or} \,\, 1, \,\,\, \forall \,\, (i,j). \label{zero1}
\end{eqnarray}
\noindent This formulation, 
of course, is impractical, because we do not have a compact way of representing any
of the constraints (\ref{nonlinear}), and there are an exponential number of
them.

Putting these issues aside, we can outline an algorithm for the min-cardinality
attack problem. Our algorithm will be iterative, 
and will maintain a ``master attacker'' mixed-integer program
which will be a {\em relaxation} of
(\ref{mincard})-(\ref{zero1}) -- i.e. it will have objective (\ref{mincard}) 
but weaker constraints than (\ref{nonlinear}).  Initially, the
master attacker MIP will include no variables other than the $z$ variables,
and no constraints other than (\ref{zero1}).  The algorithm proceeds as follows.\\

\begin{center}
  \textsc{{\bf Basic algorithm for min-cardinality attack problem}}\vspace*{5pt}\\
  \fbox{
    \begin{minipage}{0.9\linewidth}
      \hspace*{.1in} {\bf Iterate:}\\ \\
      \hspace*{.1in} {\bf 1. Attacker:} Solve master attacker MIP and let $z^*$ be its solution.\\ \\
      \hspace*{.1in} {\bf 2. Controller:} Search for a set $\cC$ of generators such that $t_{\cC}(z^*) \, \le \, 1$.\\ \\
\hspace*{.2in} {\bf (2.a)} If no such set $\cC$ exists, {\bf EXIT}:\\ \hspace*{.6in} $\sum_{ij} z^*_{ij}$ is  the minimum cardinality of a successful attack.\\ \\
\hspace*{.2in} {\bf (2.b)} Otherwise, suppose such a set $\cC$ is found.\\ \hspace*{.6in}   Add to the master attacker MIP a system of
valid inequalities that cuts off $z^*$.\\  \hspace*{.6in}   Go to {\bf 1.}

\end{minipage}
  }
\end{center}

\noindent As discussed above, the search in 
Step 2 can be implemented by solving a 
mixed integer program.  Since in 2.b we add valid inequalities to the master,
then inductively we always have a relaxation of (\ref{mincard})-(\ref{zero1})
and thus the value of the master at any execution of step 1, i.e. the
value $\sum_{ij} z^*_{ij}$, is a lower bound on the cardinality of any
successful attack.  Thus the exit condition in step 2.a is correct, since
it proves that the attack implied by $z^*$ is successful.

The implementation of Case 2.b, on the other hand, requires some care.  
Assuming we are in case 2.b, we have that $t_{\cC}  (z^*) \le  1$,
and certainly the linear program  \bmath{K_{\cC} (\cA)} is feasible.  
The optimal dual solution would therefore (apparently) furnish a Benders cut that cuts off
$z^*$. However this would be incorrect since the structure of constraints 
(\ref{2a})-(\ref{2aq}) depends on $z^*$ itself.

Instead, we need to proceed as in two-stage stochastic programming with 
recourse, where the $z$ variables play the role as ``first-stage'' variables
{\em and} also appear in the second-stage problem (the subproblem); solutions
to the dual of the second-stage problem can then be used to generate cuts
to add to the master problem.  Toward this goal,
we proceed as follows, given $\cC$ and $z^*$:
\begin{itemize}
\item[B.1] Write the {\em dual} of $K_{\cC}(\emptyset)$.
\item[B.2] As is standard in interdiction-type problems 
(see \cite{limsmith}, \cite{lind}, \cite{dnralphs}, \cite{alvarez}),
 the dual is then ``combinatorialized'' by adding the $z$ variables and
additional constraints. For example, if 
$\beta_{ij}$ indicates the dual of constraint (\ref{2a}), then we add, to the 
dual of $K_{\cC} (\emptyset)$,
inequalities of the form
$$ \beta_{ij} - M^1_{ij} z_{ij} \le M^1_{ij}, \,\,\, -\beta_{ij} - M^1_{ij} z_{ij} \le M^1_{ij},$$
\noindent for an appropriate constant $M^1_{ij} > 0$.  We proceed similarly
with constraint (\ref{2ap}), obtaining the ``combinatorial dual''.  This
combinatorial dual is the functional equivalent of the second-stage problem
in stochastic programming.
\item[B.3] Fix the $z_{ij}$ variables in the combinatorial dual to $z^*$; this
yields a problem that is equivalent to 
$K_{\cC} (z^*)$ and (using $v$ to represent the variables in aggregate form) has the general structure
\begin{eqnarray}
t_{\cC}  (z^*) \, = \, && \max \,\,\, c^T v \nonumber \\
&& P v \, \le \, b \, + \, Q z^*. \label{combdualfixed}
\end{eqnarray}
\noindent Here $P$ and $Q$ are matrices, and $b$ is
a vector, of appropriate dimensions; and we have a maximization problem
since the $K_{\cC}()$ are minimization problems.  
We obtain a cut of the form 
$$ \bar \alpha^T(b + Q z) \ge 1 + \epsilon$$ 
\noindent where $\epsilon > 0$ is a small constant and
 $\bar \alpha$ is the vector of 
optimal dual variables to (\ref{combdualfixed}). Since by assumption 
$t_{\cC}  (z^*) \le 1$ this inequality cuts off $z^*$.
\end{itemize}

\noindent Note the use of the tolerance $\epsilon$. The use of
this parameter gives {\em more} power to the controller, i.e. ``borderline''
attacks are not considered successful.  In a strict sense, therefore, we
are not solving the optimization problem to exact precision; nevertheless
in practice we expect our relaxation to have negligible impact so long as
$\epsilon$ is small.  A deeper issue here is how to interpret 
truly borderline attacks that are successful according to our strict model
(and which our use of $\epsilon$ disallows); we expect 
that in practice such attacks would be ambiguous and that the approximations
incurred in modeling power flows, estimating demands levels, 
handling fluctuating demand levels, and so on,
not to mention the numerical sensitivity of the integer and linear solvers
being used, 
would have a far more significant impact on  precision.

\subsubsection{Discussion of the initial formulation}

In order to make the outline provided in B.1-B.3 into a formal algorithm, we need to 
specify the constants $M^1_{ij}$.  As is well-known, the folklore of integer programming 
dictates that the $M^1_{ij}$ should be chosen small
to improve the {\em tightness} of the linear programming relaxation of the master
problem, that is to say, how close an approximation to the integer program
is provided by the LP relaxation.\\

While this is certainly true, we have found that, additionally, popular
optimization packages show significant numerical instability when solving
power flow {\em linear} programs.  In fact, in our experience 
it is primarily this behavior that 
mandates that the $M^1_{ij}$ should be kept as small as possible. In particular
we do not want the $M^1_{ij}$ to grow with network since this would
lead to an nonscalable approach. \\

It turns out that our formulation $K_{\cC}(\cA)$ is not ideal toward this
goal.  A particularly thorny issue is that the attack $\cA$ may disconnect the network,
and proving ``reasonable'' upper bounds on the dual variables to (for example) constraint (\ref{2b}), 
when the network is disconnected,  does not seem possible. 
In the next section we describe a different formulation for the min-cardinality
attack problem which is much better
in this regard.  Our eventual algorithm will apply
steps B.1 - B.3 to this improved formulation, while the rest of our basic
algorithmic methodology as described above will remain unchanged.

\subsection{A better mixed-integer programming formulation}\label{mip2}
As before, let $\cA$ be an attack and $\cC$ a (given) configuration of 
generators.  Let  $y^{\cC} \in R^{\cG}$ be the indicator vector 
for $\cC$, i.e. 
$y^{\cC}_i = 1$ if $i \in \cC$ and $y^{\cC}_i = 0$ otherwise. The 
linear programming
formulation presented below is similar to 
(\ref{mintobj})-(\ref{2g}) except that constraint (\ref{3fake}) is used
instead of (\ref{2aq}), and we use the indicator variables $y^{\cC}$ in
equation (\ref{powercontrol}) instead of equation (\ref{2f}). Other
differences are explained after the formulation; to 
the left of each constraint we have indicated corresponding dual variables.

\begin{eqnarray}
\mbox{\bmath{K^*_{\cC} (\cA):}} \hspace{.6in} t^*_{\cC}(z^\cA) && \doteq \ \mbox{ \ min} \ t \label{mintobj2}\\
\mbox{Subject to:} \ \ \ \ \ \ \  && \nonumber \\
\mbox{\bmath{(\alpha^{\cC}_{i})}} \ \ \ \ \ \ & & \sum_{(i,j)\in \delta^{+}(i)} f_{ij} - \sum_{(j,i)\in \delta^{-}(i)} f_{ji}  = \left \{
\begin{array}{lll}
~~P_i & {i \in \cG} &\\
-D_i & {i \in \cD} &  \\
~~~0   & \mbox{otherwise} &
\end{array} \right. \label{3b}\\
&& \nonumber \\
\mbox{\bmath{(\beta^{\cC}_{ij})}}  \ \ \ \ \ \ & &  \theta_i - \theta_j - x_{ij}f_{ij} = 0  \ \ \ \forall \ (i,j) \notin \cA \label{3a}\\
&& \nonumber \\
\mbox{\bmath{(p^{\cC}_{ij}, q^{\cC}_{ij})}} \ \ \ \ \ \  & &  u_{ij}~ t \ - \ |f_{ij}|  \ \ \ \ge \ \ \ 0, \ \ \ \forall \ (i,j) \notin \cA \label{3ap}\\
&& \nonumber \\
\mbox{\bmath{(\omega_{ij}^{\cC +}, \omega_{ij}^{\cC -})}} \ \ \ \ \ \  & &  t \ - \ |f_{ij}|  \ \ \ \ge \ \ \ 1, \ \ \ \forall \ (i,j) \in \cA \label{3fake}\\
&& \nonumber \\
\mbox{\bmath{(\gamma^{\cC+}_i, \gamma^{\cC +}_i)}}  \ \ \ \ \ \ &&   P_i^{min} y^{\cC}_i \ \leq P_i \ \leq P_i^{max} y^{\cC}_i \ \ \ \forall i \in \cG \label{3f} \label{powercontrol} \\
&& \nonumber \\
\mbox{\bmath{(\mu^{\cC})}} \ \ \ \ \ \  && \sum_{j \in \cD} D_j ~ \geq ~ T^{min} \left( \sum_{j \in \cD} D^{nom}_j \right), \label{3pg} \\
&& \nonumber \\
\mbox{\bmath{(\Delta^{\cC}_j)}} \ \ \ \ \ \  && D_j ~ \leq D_j^{nom} \ \ \ \forall j \in \cD \label{3g} \\
&& \nonumber \\
&& P \ge 0, \ \ \ D \ge 0. \label{nonneg3}
\end{eqnarray}

\noindent 
Note that (\ref{3ap}) can be represented as two linear inequalities, and the same
applies to (\ref{3fake}).  

Also, we do not force  $f_{ij} = 0$ for $(i,j) \in \cA$.  
Thus, the controller
has significantly more power than in $K_{\cC} (\cA)$.  However, because of constraint
(\ref{3fake}), we have $t^*_{\cC}(z^\cA) > 1$ as soon as any of the arcs in $\cA$ actually
carries flow.  We can summarize these remarks as follows:
\begin{RE}\label{inf2def} $\cA$ defeats $\cC$ if and only if
$t^*_{\cC}(z^\cA) > 1$.
\end{RE}

\noindent The above formulation depends on 
$\cC$ only through constraint (\ref{powercontrol}). 
Using appropriate matrices $\bar A_f, \bar A_{\theta}, \bar A_{P}, \bar A_{D}, \bar A_{t}$, and vector $\hat b$, (\ref{mintobj2})-(\ref{nonneg3}) can be abbreviated as
\begin{eqnarray*}
\mbox{\bmath{K^*_{\cC} (\cA):}} \hspace{.6in} t^*_{\cC}(z^\cA) && \doteq \ \mbox{ \ min} \ t \\
\mbox{Subject to:} \ \ \ \ \ \ \  && \nonumber \\
&& \bar A_f f \, + \bar A_{\theta} \theta \, + \bar A_P P \, + \bar A_D D \, + \bar A_t t \, \, \ge \, \, \bar b \\
 &&   P_i^{min} y^{\cC}_i \ \leq P_i \ \leq P_i^{max} y^{\cC}_i, \ \ \ \forall i \in \cG 
\end{eqnarray*}

\noindent Allowing the $y$ quantities to become $0/1$ variables, we obtain the problem
\begin{eqnarray}
\hspace{.6in} t^*(z^\cA) && \doteq \ \mbox{ \ min} \ t \label{controller-1}\\
\mbox{Subject to:} \ \ \ \ \ \ \  && \nonumber \\
&& \bar A_f f \, + \bar A_{\theta} \theta \, + \bar A_P P \, + \bar A_D D \, + \bar A_t t \, \, \ge \, \, \bar b \\
 &&   P_i^{min} y_i \ \leq P_i \ \leq P_i^{max} y_i, \ \ \ \forall i \in \cG \\
&& y_i \,\, = \,\, 0 \,\, \mbox{or} \,\, 1, \ \ \ \forall i \in \cG.  \label{controller-4}
\end{eqnarray}
\noindent This is the {\em controller's problem}: we have that $t^*(z^\cA) \le 1$ if and only if there exists some configuration
of the generators that defeats $\cA$. \\

\noindent However, for the
purposes of this section, we will assume $\cC$ is given and that the $y^{\cC}$ are
constants.  We can now write the dual of $K^*_{\cC} (\cA)$, suppressing the index $\cC$ from the variables, for clarity.

\begin{eqnarray}
\mathbf{A_{\cC}(\cA)}:  &\mbox{max}& \sum_{i \in cG} y^{\cC}_i P_i^{min} \gamma_i^{-} - \sum_{i \in \cG} y^{\cC}_i P_i^{max} \gamma_i^{+} - \sum_{j \in \cD} D_j^{nom} \Delta_j + \sum_{j \in \cD} D^{nom}_j \,\mu_j  + \sum_{(i,j) \in E} (\omega_{ij}^+ + \omega_{ij}^-) \nonumber 
\end{eqnarray}
\begin{eqnarray}
\ \ \ \ \ \ \mbox{Subject to:}  \ \ \ \  && \nonumber \\
&(f_{ij})& \alpha_i - \alpha_j - x_{ij}\beta_{ij} - p_{ij} + q_{ij} + \omega_{ij}^{+} - \omega_{ij}^{-} \,\, = \,\, 0  \ \ \forall (i,j) \in E \label{flowdual} \\
&(\theta_{i})& \sum_{(i,j) \in \delta^{+}(i)} \beta_{ij} - \sum_{(j,i) \in \delta^{-}(i)} \beta_{ji} \,\, = \,\, 0 \ \ \forall i \in V \label{thetadual} \\ 
&(t)& \sum_{(i,j) \in E} u_{ij} (p_{ij} +q_{ij}) + \sum_{(i,j) \in E} (\omega_{ij}^{+} + \omega_{ij}^{-})  \,\,\, \leq \,\,\, 1  \label{metric} \\
&(P_i)& -\alpha_i - \gamma_i^{-} + \gamma_i^{+} \,\, = \,\, 0 \ \ \forall i \in \cG \\
&(D_j)& \alpha_j  + \mu - \Delta_j \,\, \leq 0 \,\, \ \ \forall j \in \cD \\
&(\xi^+_{ij}, \xi^-_{ij})& x_{ij}^{1/2} |\beta_{ij}|  \,\, \leq \,\, \bM (1 - z^{\cA}_{ij})\ \ \forall (i,j) \in E  \label{xi}\\
&(\varrho_{ij})& p_{ij} + q_{ij} \,\, \leq \,\, \frac{1}{u_{ij}}(1 - z^{\cA}_{ij}) \ \ \forall (i,j) \in E \label{m2} \\
&(\eta_{ij})& \omega_{ij}^{+} + \omega_{ij}^{-} \,\, \leq \,\, z^{\cA}_{ij}  \ \ \forall (i,j) \in E \label{omegavub} \\
&  & \omega_{ij}^{+} \geq 0,\ \ \omega_{ij}^{-} \geq 0, \ \ p_{ij} \,\, \geq \,\, 0, \ \ q_{ij} \geq 0 \ \ \forall (i,j) \in E  \nonumber \\
&  & \gamma^+_i, \gamma^-_i \,\, \geq \,\, 0 \ \ \forall  i \in \cG \nonumber \\
&  & \Delta_j \geq 0 \ \ \forall j \in \cD  \nonumber \\
&  & \mu \geq 0 \nonumber \\
&  & \delta_{ij}, \ \beta_{ij} \ \mbox{free} \ \ \forall (i,j) \in E  \nonumber \\
&  & \alpha_{i} \ \mbox{free} \ \ \forall i \in V. \nonumber 
\end{eqnarray}
\noindent As before, for each constraint we indicate the corresponding dual variable. In (\ref{xi}), $\bM$ is an appropriately chosen constant (we will 
provide a precise value for it below).  Note
that we are scaling $\beta_{ij}$ by $x_{ij}^{1/2}$ -- this is allowable since
$x_{ij}^{1/2} > 0$; the reason for this scaling will become clear below.

\noindent Denoting
$$ \psi^{\cC} \ := \ (\alpha^{\cC}, \beta^{\cC}, p^{\cC}, q^{\cC}, \omega^{\cC +}, \omega^{\cC -}, \gamma^{\cC -}, \gamma^{\cC +}, \mu^{\cC}, \Delta^{\cC})$$ 
\noindent we have that $A_{\cC}(\cA)$ can be rewritten
as:
\begin{eqnarray}
&& \max \left\{ \, w_{\cC}^T \, \psi^{\cC} \, : \, A \psi^{\cC} \, \le \, b \, + \, B \left(1 - z^{\cA} \right) \, \right\} 
\end{eqnarray}
\noindent where $A$, $B$, $w_{\cC}$ and $b$ are appropriate matrices and vectors.  Consequently, we can now
rewrite the min-cardinality attack problem:
 \begin{eqnarray}
\mbox{ \ min} \,\, \sum_{(i,j)} z_{ij} && \label{mincard2} \\
\mbox{Subject to:}  \ \ \ \ \ \ \ \ \ \ \ \ \ \ \ \ \  t^{\cC} & \ge & 1 + \epsilon, \,\,\,\,\,\, \forall \,\, \cC \subseteq \cG \label{nonlinear0}\\
w_{\cC}^T \, \psi^{\cC} \, - \, t^{\cC} & \ge & 0, \,\,\,\,\,\, \forall \,\, \cC \subseteq \cG, \label{nonlinear2} \\
A \psi^{\cC} \, + \, B z & \le & b \, + \, B \,\,\,\,\,\, \forall \,\, \cC \subseteq \cG, \label{nonlinear3} \\
z_{ij} & = &  0 \,\, \mbox{or} \,\, 1, \,\,\, \forall \,\, (i,j). \label{zero2}
\end{eqnarray}
\noindent This formulation, of course, is exponentially large.  An alternative
is to use Benders cuts --  having solved the linear program $A_{\cC}(\cA)$,  let 
$(\bar f, \bar \theta, \bar t, \bar P, \bar D, \bar \xi^+, \bar \xi^-, 
\bar \varrho, \bar \eta)$ be optimal dual variables.  Then the resulting 
Benders cut is
\begin{eqnarray}
&& t^{\cC} + \sum_{(i,j) \in E} ( (\bar \xi^{+}_{ij} + \bar \xi^{-}_{ij}) \bM (1- z_{ij})) +  \sum_{(i,j) \in E} ( \frac{1}{u_{ij}} \bar \varrho_{ij}(1 - z_{ij})) + \sum_{(i,j) \in E} \bar \eta_{ij}z_{ij} \ge 1 + \epsilon, \label{benders-cut}
\end{eqnarray}
\noindent We can now update our algorithmic template for the min-cardinality problem.

\begin{center}
  \textsc{ {\bf Updated algorithm for min-cardinality attack problem}}\vspace*{5pt}\\
  \fbox{
    \begin{minipage}{0.9\linewidth}
      \hspace*{.1in} {\bf Iterate:}\\ \\
      \hspace*{.1in} {\bf 1. Attacker:} Solve master attacker MIP, obtaining attack $\cA$.\\ \\
      \hspace*{.1in} {\bf 2. Controller:} Solve the controller's problem (\ref{controller-1})-(\ref{controller-4}) to search for a  set \\ \hspace*{.4in}$\cC$ of generators such that $t^*_{\cC}(z^\cA) \, \le \, 1$.\\ \\
\hspace*{.2in} {\bf (2.a)} If no such set $\cC$ exists, {\bf EXIT}: \\ \hspace*{.6in}  $\cA$ is  a minimum cardinality successful attack.\\ \\
\hspace*{.2in} {\bf (2.b)} Otherwise, suppose such a set $\cC$ is found.  Then\\
\hspace*{.7in}   {\bf (2.b.1)} Add to the master the Benders' cut (\ref{benders-cut}), and, optionally\\ 
\hspace*{.7in}   {\bf (2.b.2)} Add to the master the entire system (\ref{nonlinear0})-(\ref{nonlinear3}),\\  
\hspace*{.6in}   {\bf Go to 1.}

\end{minipage}
  }
\end{center}

\noindent Clearly, option (2.b.2) can only be exercised sparingly.  Below we will 
discuss how we 
choose, in our implementation, between (2.b.1) and (2.b.2).  We will also
describe how to (significantly) strengthen the straightforward Benders cut (\ref{benders-cut}). One point to note is that the cuts (or systems)
arising from different configurations $\cC$ reinforce one another.  \\

\noindent At each iteration of the algorithm, the master attacker MIP becomes a tighter
relaxation for the min-cardinality problem, and thus its solution in step 1 provides a 
lower bound for the problem.  Thus, if in a certain execution of step 2 we certify that 
$t^*_{\cC}(z^\cA) > 1$ for every configuration $\cC$, we have solved the min-cardinality
problem to optimality.  \\

What we have above is a complete outline of our algorithm.  In order to make the algorithm
effective we need to further sharpen the approach.  In particular, we need 
set the constant $\bM$ to as small a value as possible, and we need to develop stronger
inequalities than the basic Benders' cuts.

\subsubsection{Setting  $\bM$}
\noindent In this section we show how to set for $\bM$ a value that does not grow with network size.

\begin{LE}\label{bigM}
In formulation $\mathbf{A_{\cC}(\cA)}$, a valid choice for $\bM$ is
\begin{eqnarray}
&& \bM \, = \,   \max_{(i,j) \in E} \left\{\frac{1}{\sqrt{x_{ij}}\, u_{ij}}\right\}. \label{beta-M}
\end{eqnarray}
\end{LE}
\noindent {\em Proof.} Given an attack $\cA$, consider a connected component
$K$ of $\cN - \cA$.  For any arc $(i, j)$ with both ends in $K$, 
$\omega_{ij}^{+} + \omega_{ij}^{-} = 0$ by (\ref{omegavub}).  Hence we can
rewrite constraints (\ref{flowdual})-(\ref{thetadual}) over all arcs with 
both ends in $K$ as follows:
\begin{eqnarray}
N^T_K \alpha_K \, - \, X_K \beta_K & = & p_K - q_K, \label{sys1} \\
N_K \beta_K & = & 0.  \label{sys2}
\end{eqnarray}
Here, $N_K$ is the node arc incidence matrix of $K$, $\alpha_K, \beta_K, p_K, q_K$ are the restrictions of $\alpha, \beta, p, q$ to $K$, and $X_K$ is 
the diagonal matrix $diag\{x_{ij} \, : \, (i,j) \in K \}$. From this 
system we obtain
\begin{eqnarray}
N_K X^{-1}_K N_K \alpha_K & = & N_K X^{-1}_K (p_K - q_K).  \label{sys3}
\end{eqnarray}
\noindent The matrix $N_K X^{-1}_K N_K$ has one-dimensional null space and
thus we have one degree of freedom in choosing $\alpha_K$.  Thus, to
solve (\ref{sys3}), we can remove from $N_K$ an arbitrary row, obtaining
$\tilde N_K$, and  remove the same row from $\alpha_K$, obtaining 
$\tilde \alpha_K$.  Thus, (\ref{sys3}) is equivalent to:
\begin{eqnarray}
\tilde N_K X^{-1}_K \tilde N_K \tilde \alpha_K & = & \tilde N_K X^{-1}_K (p_K - q_K),  \label{sys4}
\end{eqnarray}
\noindent The matrix $\tilde N_K X^{-1}_K \tilde N_K$ and thus (\ref{sys4})
has a unique solution (given $p_K - q_K$); we complete this to a solution to 
(\ref{sys3}) by setting to zero the entry of $\alpha_K$ that was removed.
Moreover,
\begin{eqnarray}
&& X_K^{-1/2} N^T_K \alpha_K \, = \, X_K^{-1/2} \tilde N^T_K \tilde \alpha_K \, = \,
X_K^{-1/2} \tilde N_K^T (\tilde N_K X_K^{-1} \tilde N_K^T)^{-1} \tilde N_K X_K^{-1} (p_K - q_K). \label{sys5}
\end{eqnarray}
\noindent The matrix
$$ H \, := \, X_K^{-1/2} \tilde N_K^T ~ (\tilde N_K X_K^{-1} \tilde N_K^T)^{-1} ~ \tilde N_K X_K^{-1/2}$$
\noindent is symmetric and idempotent, e.g. $H H^T = I$.  Thus, from (\ref{sys5})
we get
\begin{eqnarray}
\| X_K^{-1/2} N^T_K \alpha_K \|_2 & \le & \| H \|_2 ~ \| X_K^{-1/2} (p_K - q_K) \|_2 \, \le \, \| X_K^{-1/2} (p_K - q_K) \|_2,  \label{longone}
\end{eqnarray}
\noindent where the last inequality follows from the idempotent attribute. 
Because of constraints (\ref{metric}), (\ref{m2}) and (\ref{omegavub}), we can
see that the square of the right-hand side of (\ref{longone}) is upper-bounded by the
value of the convex maximization problem,
\begin{eqnarray}
&\mbox{max}& \sum_{(i,j) \in E} x_{ij}^{-1} (p_{ij}-q_{ij})^2 \\
&\mbox{s.t.}& \sum_{(i,j) \in E} u_{ij} (p_{ij} + q_{ij}) \leq 1 \\
&  & p_{ij} \geq 0, \ q_{ij} \geq 0,
\end{eqnarray}
\noindent which equals
$$\max_{(i,j) \in E }\left\{\frac{1}{x_{ij}u_{ij}^2}\right\}. $$
\QED

\subsubsection{Tightening the formulation}

In this section we describe a family
of inequalities that are valid for the attacker problem.  These cuts seek to capture the interplay between the 
flow conservation equations and Ohm's law.  First we present a technical result.

\begin{LE} \label{l-norms}Let $Q$ be a matrix with $r$ rows with rank $r$, 
and let $A = Q^{T}(Q Q^{T})^{-1} Q  \in \cR^{r \times r}$. Let $ B := I - A$. 
Then for any $p \in \cR^r$ we have
\begin {eqnarray}
&\|p\|_2^{2}& = \|Ap\|_2^{2} + \|Bp\|_2^{2} \label{l-norms1}\\
& \|p\|_1& \geq |(Ap)_j| + |(Bp)_j| \ \ \forall j = 1 \ldots r \label{l-norms2}
\end{eqnarray} \end{LE}

\noindent {\em Proof.} $A$ and $B$ are symmetric and idempotent, i.e., $A^{2} = A$, $B^{2} = B$, and any $p \in \cR^r$ can be written as $ p = Ap + Bp $. 
Multiplying this equation by $p$ and using the fact that $A$ and  $B$ are symmetric and idempotent we get (\ref{l-norms1}):
\begin {eqnarray}
p^{T}p &=& p^{T}Ap + p^{T}Bp \\
       &=& p^{T}A^{2}p + p^{T}B^{2}p \\
\|p\|_2^{2} &=& \|Ap\|_2^{2} + \|Bp\|_2^{2} \label{givemeabreak}
\end{eqnarray}

\noindent We also have $A^{T}B = A (I - A) = A - A^{2} = 0$, so $y^{T}A^{T}By = 0 $ for any $y \in \cR^r$. Thus, if we rename $Ap = x$ and $Bp = y$, 
then the following holds: 
$p = x + y, \ x^{T}y = 0, \ \|p\|_2^{2} = \|x\|_2^{2} + \|y\|_2^{2}$. 

Let $1 \le j \le r$.  We have
$$ \|p\|_2^{2} - (|x_j| + |y_j|)^{2} =  \|x\|_2^{2} + \|y\|_2^{2} - (|x_j| + |y_j|)^{2} = \sum_{i, i\neq j}x_i^{2} + \sum_{i, i\neq j}y_i^{2} - 2 |x_j y_j| $$

\noindent where the first equality follows from (\ref{givemeabreak}). Since $x^{T}y = 0 $, we have $|x_j y_j| = |\sum_{i, i\neq j} x_i y_i |$. Hence,
\begin {eqnarray}
\sum_{i, i\neq j}x_i^{2} + \sum_{i, i\neq j}y_i^{2} - 2 |x_j y_j| &=& \sum_{i, i\neq j}x_i^{2} + \sum_{i, i\neq j}y_i^{2} - 2 \left |\sum_{i, i\neq j} x_i y_i \right |\\
&\geq& \sum_{i, i\neq j}x_i^{2} + \sum_{i, i\neq j}y_i^{2} - 2 \sum_{i, i\neq j} | x_i y_i | \\
&=& \sum_{i, i\neq j}(|x_i| - |y_i|)^{2} \\
&\geq& 0
\end{eqnarray}
\noindent So we have $\|p\|_2^{2} - (|x_j| + |y_j|)^{2} \geq 0$, 
which implies $ \|p\|_1 \geq \|p\|_2 \geq (|x_j| + |y_j|) \ \ \forall j = 1 \ldots r$. \QED \\

\noindent As a consequence of this result we now have:

\begin{LE} \label{l-valid} Given configuration $\cC$, the following inequalities are valid for 
system (\ref{nonlinear3})-(\ref{zero2}) for each $(i,j) \, \in \, E$:
\begin {eqnarray}
&&x_{ij}^{-\frac{1}{2}} |\alpha^{\cC}_i - \alpha^{\cC}_j| + x_{ij}^{\frac{1}{2}} |\beta^{\cC}_{ij}|  \leq x_{ij}^{-\frac{1}{2}} w^{\cC}_{ij}  +  \bM (1 - z_{ij})
\label{l-valid1}\\
&&x_{ij}^{-\frac{1}{2}} |\alpha^{\cC}_i - \alpha^{\cC}_j| + x_{ij}^{\frac{1}{2}} |\beta^{\cC}_{ij}|  \leq \sum_{(k,l)}x_{kl}^{-\frac{1}{2}}( p^{\cC}_{kl} + q^{\cC}_{kl} ) +  w^{\cC}_{ij}
\label{l-valid2}
\end{eqnarray}
where $\bM := max_{(k,l) \in E}\{\frac{1}{\sqrt{x_{kl}}u_{kl}}\}$ as before.
\end{LE}
\noindent {\em Proof.} Suppose first that $z_{ij} = 0$.  Let $K$ be the component containing $(i,j)$ after
the attack. Then by (\ref{sys5}) and (\ref{sys1}),
\begin{eqnarray}
 X^{-1/2} N_K^{T}\alpha^{\cC} & = &A X^{-1/2}(p^{\cC} - q^{\cC}), \\
 X^{1/2} \beta^{\cC} & = & (I-A) X^{-1/2}(p^{\cC} - q^{\cC}),
\end{eqnarray}
\noindent where $A = X^{-1/2}\tilde{N_K}^{T}(\tilde{N_K}X^{-1}\tilde{N_K}^{T})^{-1}\tilde{N_K}X^{-1/2}$.
Thus, we have 
\begin{eqnarray}
x_{ij}^{-1/2} |\alpha^{\cC}_i - \alpha^{\cC}_j| + x_{ij}^{1/2} |\beta^{\cC}_{ij}| & \leq & \sum_{(k,l)} x_{kl}^{-1/2}( p^{\cC}_{kl} + q^{\cC}_{kl}) \,\, \leq \,\, \bM 
\end{eqnarray}
\noindent where the first inequality follows from (\ref{l-norms2}) proved in Lemma \ref{l-norms}, 
and the second bound is obtained as in the proof of Lemma \ref{bigM}. \\

\noindent Suppose now that $z_{ij} = 1$. Here we have 
$|\alpha^{\cC}_i - \alpha^{\cC}_j| \leq \omega^{\cC}_{ij}$, by (\ref{flowdual}), (\ref{m2}), (\ref{xi}).
Using these (\ref{l-valid1})-(\ref{l-valid2}) can be easily shown. \QED\\

Inequalities (\ref{l-valid1})-(\ref{l-valid2}) strengthen system (\ref{nonlinear3})-(\ref{zero2}); 
when case step (2.b.2) of the min-cardinality algorithm is applied then (\ref{l-valid1}), (\ref{l-valid2})
will become part of the master problem.  If case (2.b.1) is applied, then 
the vector $ \psi^{\cC} = (\alpha^{\cC}, \beta^{\cC}, p^{\cC}, q^{\cC}, \omega^{\cC +}, \omega^{\cC -}, \gamma^{\cC -}, \gamma^{\cC +}, \mu^{\cC}, \Delta^{\cC})$ is expanded by adding two new dual variables per arc
$(i, j)$.

\subsubsection{Strengthening the Benders cuts}\label{strengthening}
Typically, the standard Benders cuts (\ref{benders-cut}) prove weak.  One manifestation
of this fact is that in early iterations of our algorithm for the min-cardinality attack problem, the attacks produced in Step 1 will tend to be ``weak'' and, in particular, of
very small cardinality.  Here we discuss
two routines that yield substantially stronger inequalities, still in the Benders mode.  

In Step 2 of the algorithm, given an attack $\cA$, we discover a generator configuration $\cC$ that defeats
$\cA$, and from this configuration a cut is obtained.  However, it is not simply the
configuration that defeats $\cA$, but, rather, a vector of power flows.  If we could
somehow obtain a ``stronger'' vector of power flows, the resulting cut should 
prove tighter.  
To put it differently, a vector of power flows that are in some sense ``minimal'' 
might also defeat other attacks $\cA'$ that are ``stronger`` than $\cA$; in other words, 
they should produce
stronger inequalities.   \\

We implement this rough idea in two different ways.  Consider Step 2 of the min-cardinality
attack algorithm, and suppose case (2.b) takes place.  We execute steps I and II below,  where
in each case $\cA^*$ is initialized as $E - \cA$, and $f^*$ is initialized as the 
power flow that defeated $\cA$:
\begin{itemize}
\item [{\bf (I)}] First, we add the Benders' cut (\ref{benders-cut}). \\ 
  \noindent Also, initializing $\cB = \cA$,  we run the following step,  for $k = 1, 2, \ldots, |E - \cA|$:
\begin{itemize}
  \item [{\bf (I.0)}] Let $(i_k, j_k) \, = \, \argmin \left\{ |f^*_{ij}| \, : \, (i,j) \, \in \, \cA^* \right\}$.
  \item [{\bf (I.1)}] If the attack $\cB \cup (i_k, j_k)$ is {\em not} successful, then reset
    $\cB \, \leftarrow \, \cB \cup (i_k, j_k)$, and update $f^*$ to the power flow that
defeats the (new) attack $\cB$.
  \item [{\bf (I.2)}]  Reset $\cA^* \leftarrow\cA^* - (i_k, j_k)$.
\end{itemize}
\noindent At the end of the loop, we have an attack $\cB$ which is not successful, i.e.
$\cB$ is defeated by some configuration $\cC'$.  If $\cB = \cA$ we do nothing.  
Otherwise, we add to the master problem the Benders cut arising from $\cB$ and $\cC'$.
\item [{\bf (II)}] Set $\cF = \emptyset$ and $\cC' = \cC$.  We run the following step,  for $k = 1, 2, \ldots, |E - \cA|$:
\begin{itemize}
  \item [{\bf (II.0)}] Let $(i_k, j_k) \in \cA^*$ be such that its flow has minimum absolute value.
  \item[{\bf (II.1)}] Test whether $\cA$ is successful against a controller that is forced
to satisfy the condition
\begin{eqnarray}
&& f_{ij} \, = \, 0, \,\,\, \forall \,\,\, (i,j) \, \in \, \cF \cup (i_k, j_k). \label{forced}
\end{eqnarray}
\item[{\bf (II.2)}] If {\em not} successful, let $\cC'$ be the configuration that defeats
the attack,  and reset $f^*$ to the corresponding power flow that
satisfies (\ref{forced}). Reset $\cF \leftarrow \cF \cup (i_k, j_k)$, 
  \item [{\bf (II.3)}] Reset $\cA^* \leftarrow\cA^* - (i_k, j_k)$.
\end{itemize}

\end{itemize}

\noindent {\bf Comment.}  Procedure (I) produces attacks of increasing cardinality.  
At termination, if $\cA \neq \cB$, then and $\cC \neq \cC'$, 
and yet $\cB$ is still not successful.  In some
sense in this case $\cC'$ is a 'stronger' configuration than $\cC$ and the resulting Benders'
cut 'should' be tighter than the one arising from $\cC$ and $\cA$.  We say 'should' because
the previously discussed non-monotonicity property of power flow problems could mean that
$\cC'$ does not defeat $\cA$.  Nevertheless, {\em in general}, the new cut is indeed
stronger. 

In contrast with (I), procedure (II) considers a progressively weaker controller.  In
fact, because we are forcing flows to zero, but we are not voiding Ohm's equation (\ref{ohm-eq}), the power flow that defeats $\cA$ while satisfying (\ref{forced}) is a feasible power
flow for the original network. Thus, at termination of the loop,
$$ \cC' \,\,\, \mbox{defeats every attack $\cA'$ of the form $\cA' \, = \, \cA \cup \cE$ for each $\cE \subseteq \cF$}.$$ 
\noindent Thus, if $\cF \neq \emptyset$ the cut obtained in (II) should be particularly
strong.\\

One final comment on procedures (I) and (II) is that each ``test'' requires the solution
of the controller's problem (\ref{controller-1})-(\ref{controller-4}), a mixed-integer program.  In our testing, on networks with up to a few hundred arcs
and nodes, such problems can be solved in sub-second time using a commercial
integer programming solver.

\subsection{Implementation details}
Our implementation is based on the updated algorithmic outline given in Section \ref{mip2}.
In step (2.b.1) we add the Benders' cut with strengthening as in section \ref{strengthening},
so we may add two cuts.  We execute Step (2.b.2) so that the relaxation includes up to
two full systems (\ref{nonlinear0})-(\ref{nonlinear3}) at any time: when a system is added
at iteration $k$, say, it is replaced at iteration $k + 4$ by the system corresponding to the
configuration $\cC$ discovered in Step 2 of that iteration.  Because at each iteration
the cut(s) added in step (2.b.1) cut-off the current vector $z^{\cA}$, the procedure is
guaranteed to converge.

\subsection{Computational experiments with the min-cardinality model}\label{comp-model1} 
\noindent In the experiments reported in this section we used a 3.4 GHz Xeon
machine with 2 MB L2 cache and 8 GB RAM.  All experiments were run using
a single core. The LP/IP solver was Cplex v. 10.01, with default settings.   Altogether, we report on $118$ runs of our algorithm.  

\subsubsection{Data sets}

For our experiments we
used problem instances of two types; all problem instances are available for download (http://www.columbia.edu/$\sim $dano/research/pgrid/Data.zip).
\begin{itemize}
\item [(a)] Two of the IEEE ``test cases'' \cite{IEEE}: the ``57 bus'' case
($57$ nodes, $78$ arcs, $4$ generators) and the ``118 bus'' case 
($118$ nodes, $186$ arcs, $17$ generators).  
\item [(b)] Two artificial examples were also created.
One was a ``square grid'' network with $49$ nodes and $84$ arcs, $4$
generators and $14$  demand nodes. We also considered a modified 
version of this
data set with $8$ generators but equal sum $\sum_{i \in \cG} P_i^{max}$.
We point out that 
square grids frequently arise as 
difficult networks for combinatorial problems; they are sparse while at
the same time the ``squareness'' gives rise to symmetry.  We created a 
second artificial network by taking two
copies of the $49$-node network and adding a random set of 
arcs to connect the two copies;
with resistances (resp. capacities) equal to the average in the $49$-node
network plus a small random perturbation. 
This yielded a $98$- node, $204$-arc network, with $28$ demand nodes,
and we used $10$, $12$, and $15$ generator variants.

\end{itemize}
\noindent In all cases, each of the generator output lower bounds $P_i^{min}$ was set to a random fraction (but never higher than $80\%$) 
of the corresponding $P_i^{max}$.\\
 
\noindent An important consideration involves the capacities $u_{ij}$
 -- should capacities be too small, or too large, 
the problem we study tends to become quite easy (i.e. the network is
either trivially too tightly capacitated, or has very large capacity
surpluses).  For example, if a generator accounts for $20\%$ of all demand then in a tightly capacitated situation the removal of just one arc incident
with that node could constitute a successful attack for $T^{min}$ large.
For the purposes of our study, we assumed {\em constant} capacities
for the two networks in (a) and the initial network in (b); these constants
were scaled, through experiments with our algorithm, precisely to make the
problems we solve more difficult.  A topic of further research would be to analyze
the $N-k$ problem under regimes where capacities are significantly different
across arcs, possibly reflecting a condition of pre-existing stress.  
In Section \ref{expers}, which addresses experiments involving 
the second model in this paper, we consider some variations in capacities.

\subsubsection{Goals of the experiments}

The experiments focus, primarily, on the computational workload incurred by
our algorithm. First, does the running time, and, in 
particular, the number of iterations, grow very rapidly with 
network size?  Second, does the number of generators exponentially
impact performance -- does the algorithm need to enumerate a large fraction
of the generator configurations?  In general, 
what features of a problem instance adversely
affect the algorithm -- i.e., is there any particular
pattern among the more difficult cases we observe?

As noted above, previous studies (see \cite{arrgal05}, \cite{alvarez}, \cite{swb}) involving integer programming methods
applied to the $N - k$ problem have considered examples
with up to $79$ arcs (and sparse). In this paper,
 in addition to considering significantly larger examples 
(from a combinatorial standpoint) we also face the added combinatorial 
complexity caused by the generator configurations.  Potentially, therefore, our algorithm could rapidly break down -- thus, our focus on performance.

\subsubsection{Results}

Tables \ref{57bus1}-\ref{gen2} contain our results;   Tables \ref{57bus1} and \ref{118bus1} refer to the 57-bus and 118-bus case, respectively,
Table \ref{square} considers the artificial $49$-node case and
Table \ref{gen2} considers the $98$-node case. In the tables, each row corresponds to a value of the minimum throughput $T^{min}$, while
each column corresponds to an attack cardinality.  For each (row, column)
combination, the corresponding cell is labeled ``Not Enough'' when 
using any attack of 
the corresponding cardinality (or smaller) the attacker will not be able to
reduce demand below the stated throughput, while ``Success'' means that some
attack of the given cardinality (or smaller) does succeed.  Further, we
also indicate the number of iterations 
that the algorithm took in order to prove the given outcome (shown in parentheses) as well as the corresponding CPU time in seconds. Thus, for example, in Table \ref{118bus1}, the algorithm {\em proved} that
using an attack of size $3$ or smaller we cannot reduce total demand below 
$75 \%$ of the nominal value; this required $4$ iterations which overall
took $267$ seconds. At the same time, in $7$ iterations ($6516$ seconds)
the algorithm found a successful attack of cardinality $4$.

Comparing the 57- and 118-bus cases, the 
significantly higher CPU times for the second case could be explained
by the much larger number of arcs.  The larger number of generators
could also be a cause -- however,
the number of generator configurations in the second case is more than 
eight thousand
times larger than that of the first; much larger than the actual slowdown 
shown by the tables.  

\begin{table}[h]
\vskip 14 pt
\centering
\caption{{\bf \emph{Algorithm for min-cardinality problem on 57-bus test case}}} 
\begin{tabular}{|c| c| c| c| c| c|}
\hline
 \multicolumn{6}{|c|}{{\bf 57 nodes, 78 arcs, 4 generators}}\\
 \multicolumn{6}{|c|}{Entries show: (iteration count), CPU seconds,}\\
 \multicolumn{6}{|c|}{Attack status ({\bf F =} cardinality too small, {\bf S} = attack success)}\\
\hline
\hline
 \multicolumn{1}{|c|} {} & \multicolumn{5}{|c|}{{\bf Attack cardinality}}\\\hline
{\bf Min. throughput} & {\bf 2} & {\bf 3} & {\bf 4} & {\bf 5}  & {\bf 6} \\\hline
{\bf 0.75}  & (1), 2, {\bf F} & (2), 3, {\bf S} &  & & \\
\hline
{\bf 0.70}  & (1), 1, {\bf F} & (3), 7, {\bf F} & (48), 246, {\bf F} & (51), 251, {\bf S} & \\
 \hline
{\bf 0.60}  & (2), 2, {\bf F} & (3), 6, {\bf F} & (6), 21, {\bf F} & (6), 21, {\bf S}& \\
\hline
{\bf 0.50}  & (2), 2, {\bf F} & (3), 7, {\bf F} & (6), 13, {\bf F}  &(6), 13, {\bf F} & (6), 13, {\bf S}\\
\hline
{\bf 0.30} & (1), 1, {\bf F} & (2), 3, {\bf F} & (2), 3, {\bf F} & (2), 3, {\bf F} & (2), 3, {\bf F} \\
\hline
\end{tabular}
\label{57bus1}
\end{table}

\begin{table}[h]
\vskip 14 pt
\centering
\caption{{\bf \emph{Algorithm for min-cardinality problem on 118-bus test case}}} 
\begin{tabular}{|c| c| c| c|}
\hline
 \multicolumn{4}{|c|}{{\bf 118 nodes, 186 arcs, 17 generators}}\\
 \multicolumn{4}{|c|}{Entries show: (iteration count), CPU seconds,}\\
 \multicolumn{4}{|c|}{Attack status ({\bf F =} cardinality too small, {\bf S} = attack success)}\\
\hline
\hline
 \multicolumn{1}{|c|} {} & \multicolumn{3}{|c|}{{\bf Attack cardinality}}\\\hline
{\bf Min. throughput} & {\bf 2} & {\bf 3} & {\bf 4} \\\hline
{\bf 0.92}  & (4), 18, {\bf S} & &  \\
\hline
{\bf 0.90}  & (5), 180, {\bf F} & (6), 193, {\bf S} &   \\
\hline
{\bf 0.88}  & (4), 318, {\bf F} & (6), 595, {\bf S} &   \\
 \hline
{\bf 0.84}  & (2), 23, {\bf F} & (6), 528, {\bf F} & (148), 6562, {\bf S}  \\
\hline
{\bf 0.80}  & (2), 18, {\bf F} & (5), 394, {\bf F} & (7), 7755, {\bf F}  \\
\hline
{\bf 0.75} & (2), 14, {\bf F} & (4), 267, {\bf F} & (7), 6516, {\bf F} \\
\hline
\end{tabular}
\label{118bus1}
\end{table}

\begin{table}[tbh]
\vskip 14 pt
\centering
\caption{{\bf \emph{Algorithm for min-cardinality problem on small network}}}
\begin{tabular}{|c| c| c| c| c|}
\hline
 \multicolumn{5}{|c|}{{\bf 49 nodes, 84 arcs}}\\
 \multicolumn{5}{|c|}{Entries show: (iteration count), CPU seconds,}\\
 \multicolumn{5}{|c|}{Attack status ({\bf F =} cardinality too small, {\bf S} = attack success)}\\
\hline \hline
 \multicolumn{5}{|c|}{{\bf 4 generators}}\\
\hline 
 \multicolumn{1}{|c|} {} & \multicolumn{4}{|c|}{{\bf Attack cardinality}}\\
\hline
\textbf{Min. throughput}  & {\bf 2} & {\bf 3} & {\bf 4} & {\bf 5} \\
\hline 
{\bf 0.84} & (4), 129, {\bf F} & (4), 129, {\bf S} &  & \\
\hline
{\bf 0.82}  & (4), 364, {\bf F} & (35), 1478, {\bf F}& (36), 1484, {\bf S} & \\
\hline 
{\bf 0.78}  & (4), 442, {\bf F} & (4), 442, {\bf F} & (26), 746, {\bf S} & \\
\hline 
{\bf 0.74}  & (4), 31, {\bf F} & (11), 242, {\bf F} & (168), 4923, {\bf F} & (168), 4923, {\bf S}\\
\hline 
{\bf 0.70}  & (3), 31, {\bf F} & (4), 198, {\bf F} & (10), 1360, {\bf F} & (203), 3067, {\bf S}\\
\hline 
{\bf 0.62}  & (4), 86, {\bf F} & (4), 86, {\bf F} & (131), 2571, {\bf F} & (450), 34298, {\bf F}\\
\hline
\hline 
 \multicolumn{5}{|c|}{{\bf 8 generators}}\\
\hline \hline
 \multicolumn{1}{|c|} {} & \multicolumn{4}{|c|}{{\bf Attack cardinality}}\\
\hline
\textbf{Min. throughput} & {\bf 2} & {\bf 3} & {\bf 4} & {\bf 5} \\
\hline \hline
{\bf 0.90}  & (1), 13, {\bf F} & (3), 133, {\bf S} &  & \\
\hline
{\bf 0.86}  & (1), 59, {\bf F} & (5), 357, {\bf F} & (13), 1291, {\bf S} & \\
\hline
{\bf 0.84}  & (1), 48, {\bf F} & (4), 227, {\bf F} & (41), 2532, {\bf F} & (43), 2535, {\bf S}\\
\hline
{\bf 0.80}  & (1), 14, {\bf F} & (4), 210, {\bf F} & (8), 1689, {\bf F} & (50), 2926, {\bf S}\\
\hline
{\bf 0.74} & (1), 8, {\bf F} & (3), 101, {\bf F} & (10), 1658, {\bf F} & (68), 23433, {\bf F}\\
\hline
\end{tabular}
\label{square}
\end{table}

\begin{table}[tbh]
\vskip 14 pt
\centering
\caption{{\bf \emph{Algorithm for min-cardinality problem on larger network}}}
\begin{tabular}{|c| c| c| c|}
\hline
 \multicolumn{4}{|c|}{{\bf 98 nodes, 204 arcs}}\\
 \multicolumn{4}{|c|}{Entries show: (iteration count), CPU seconds,}\\
 \multicolumn{4}{|c|}{Attack status ({\bf F =} cardinality too small, {\bf S} = attack success)}\\
\hline 
\hline
 \multicolumn{4}{|c|}{{\bf 10 generators}}\\\hline
 \multicolumn{1}{|c|} {} & \multicolumn{3}{|c|}{{\bf Attack cardinality}}\\\hline
{\bf Min. throughput} & {\bf 2} & {\bf 3} & {\bf 4} \\\hline
 {\bf 0.89} & (2) 177, {\bf F} & (30) 555, {\bf S} &  \\
\hline
{\bf 0.86}  & (2), 195, {\bf F} & (12), 5150, {\bf F} & (14), 5184, {\bf S} \\
 \hline
{\bf 0.84}  & (2), 152, {\bf F} & (11), 7204, {\bf F} & (35), 223224, {\bf F} \\
\hline
{\bf 0.82} & (2), 214, {\bf F} & (9), 11458, {\bf F} & (16), 225335, {\bf F}  \\
\hline
{\bf 0.75} & (2), 255, {\bf F} & (9), 5921, {\bf F} & (17), 151658, {\bf F} \\
\hline
{\bf 0.60}  &  & (1), 4226, {\bf F} & N/R \\
\hline 
\hline
 \multicolumn{4}{|c|}{{\bf 12 generators}}\\\hline
 \multicolumn{1}{|c|} {} & \multicolumn{3}{|c|}{{\bf Attack cardinality}}\\\hline
{\bf Min. throughput} & {\bf 2} & {\bf 3} & {\bf 4} \\\hline
{\bf 0.92}  & (2), 318, {\bf F} & (11), 7470, {\bf F} & (14), 11819, {\bf S}  \\
\hline
{\bf 0.90}  & (2), 161, {\bf F} & (11), 14220, {\bf F} & (18), 16926, {\bf S} \\
 \hline
{\bf 0.88}  & (2), 165, {\bf F} & (10), 11178, {\bf F} & (15), 284318, {\bf S} \\
\hline
{\bf 0.84}  & (2), 150, {\bf F} & (9), 4564, {\bf F} & (16), 162645, {\bf F}  \\
\hline
{\bf 0.75} & (2), 130, {\bf F} & (9), 7095, {\bf F} & (15), 93049, {\bf F} \\
\hline
\hline
 \multicolumn{4}{|c|}{{\bf 15 generators}}\\\hline
 \multicolumn{1}{|c|} {} & \multicolumn{3}{|c|}{{\bf Attack cardinality}}\\\hline
{\bf Min. throughput} & {\bf 2} & {\bf 3} & {\bf 4} \\\hline

{\bf 0.94}  & (2), 223, {\bf F} & (11), 654, {\bf S} &  \\
\hline
{\bf 0.92}  & (2), 201, {\bf F} & (11), 10895, {\bf F} & (18), 11223, {\bf S} \\
 \hline
{\bf 0.90}  & (2), 193, {\bf F} & (11), 6598, {\bf F} & (16), 206350, {\bf S} \\
\hline
{\bf 0.88}  & (2), 256, {\bf F} & (9), 15445, {\bf F} & (18), 984743, {\bf F}  \\
\hline
{\bf 0.84}  & (2), 133, {\bf F} & (9), 5565, {\bf F} & (15), 232525, {\bf F}  \\
\hline
{\bf 0.75}  & (2), 213, {\bf F} & (9), 7550, {\bf F} & (11), 100583, {\bf F} \\
\hline
\end{tabular}
\label{gen2}
\end{table}

Table \ref{square} presents experiments with our
algorithm on the $49$-node, $84$-arc network, first using $4$ and then
$8$ generators.  Not surprisingly, the network with $8$ generators proves more resilient (even though total generator capacity is the same) -- for
example, an attack of cardinality $5$ is needed to reduce throughput below
$84 \%$, whereas the same can be achieved with an attack of size $3$ in the
case of the $4$-generator network.  Also note that the running-time performance
does not significantly degrade as we move to the $8$-generator case, even
though the number of generator configurations has grown by a factor of
$16$.  Not surprisingly,
the most time-consuming cases are those where the adversary fails, since
here the algorithm must prove that this is the case (i.e. prove that no
successful attack of a given cardinality exists) while in a ``success''
case
the algorithm simply needs to find {\em some} successful attack of the right
cardinality.\\

Table \ref{gen2} describes similar tests, but now on the
$98$-node, $204$-arc network.  Note that in the $15$ generator case
there are over $30000$ generator configurations that must be examined,
at least implicitly, in order to certify that a given attack is successful. 
But, as in the case of Table \ref{square}, the number of generators
does not have an exponential impact on the overall running time.  \\

In general, therefore,  the experiments show that
the number
of generators plays a second-order role in the complexity of the algorithm;
the total number of iterations depends weakly on the total number of
generator configurations, and the primary agent behind complexity is the topological network 
structure.\\

A point
worth dwelling on is that, with a few exceptions, the running time tends to
{\em decrease}, for a given attack cardinality, once the minimum throughput
is sufficiently past the threshold where no successful attack exists.  This
can be explained as follows: as the minimum throughput decreases the controller
has more ways to defeat the attacker -- if no attack can succeed a pure
enumeration algorithm would have to enumerate all possible attacks; thus
arguably our cutting-plane approach does indeed discover useful structure
that limits enumeration (i.e., the cuts added in step (2.b.1) of our algorithm
enable us to prove an effective lower bound on the minimum attack 
cardinality needed to obtain a successful attack).

Also (consider the cases corresponding to cardinality = $4$; and we have $12$ or $15$ 
generators) 
 CPU time increases with decreasing minimum throughput so long as a 
successful attack {\em does} exist.  This can also be explained, 
as follows: in order
for the algorithm to terminate it must generate a successful attack, 
but this search becomes more difficult as the minimum throughput 
decreases (the controller has more options).  Roughly speaking, in summary, we
would expect the problem  to be ``easiest'' (for a given attack cardinality ) 
near extreme values of the
minimum demand threshold; the experiments tend to confirm this expectation.

\begin{table}[tbh]
\centering
\caption{{\bf Min-cardinality algorithm on one-configuration problem}} 
\begin{tabular}{|c|c|c|}
\hline
{\bf Min. Throughput} & {\bf Min. Attack Size} & {\bf Time (sec.)} \\
\hline
0.95 & 2 & 2 \\
0.90 & 3 & 20 \\
0.85 & 4 & 246 \\
0.80 & 5 & 463 \\
0.75 & 6 & 2158 \\
0.70 & 6 & 1757 \\
0.65 & 7 & 3736 \\
0.60 & 7 & 1345 \\
0.55 & 8 & 2343 \\
0.50 & 8 & 1328 \\
\hline
\end{tabular}
\label{onescenario}
\end{table}

\subsubsection{One-configuration problems}
For completeness, in Table \ref{onescenario} we present results where
we study {\em one-configuration} problems where the set of generators
that the controller operates are {\em fixed}.  For a given minimum demand
throughput, the table shows the minimum attack cardinality needed to defeat
the controller. Problems of this type correspond
most closely to those previously studied in the literature. Here we 
applied the mixed-integer programming formulation (\ref{mincard2})-(\ref{zero2}) restricted to
the single configuration $\cC = \cG$.  Rather than use our algorithm, we 
simply solved these problems using Cplex, with default settings.  The table
shows the CPU time needed to solve the minimum-cardinality problem 
corresponding to the minimum throughput shown in the first column.  The
point here is that our formulation (\ref{mincard2})-(\ref{zero2}) proves 
significantly effective in relation to previous methods.

Not surprisingly, the problem becomes {\em easier} as the attack cardinality increases 
-- more candidates (for optimal attack) exist.

\section{A continuous, nonlinear attack problem}\label{nlp}

In this section we study a new attack model.  Our goals are twofold:
\begin{itemize}
\item First, we want to more explicitly capture how the flow conservation equations (\ref{1b})
interact with the power-flow model (\ref{ohm-eq}) in order to produce flows in
excess of capacities. More generally, we are interested in directly 
incorporating the interaction
of the laws of physics with the graph-theoretic structure of the network into an 
algorithmic procedure.  
It is quite clear that the complexity of combinatorial
problems on power flows, such as the min-cardinality attack problem, is
primarily due to this interaction.  
\item Second, there are ways other than the outright disabling of a power line, in which
the functioning of the line could be hampered.  There is a sense (see e.g. \cite{usc})
that recent real-world blackouts were not simply the result of discrete line
failures; rather the system as a whole was already under ``stress'' when the
failures took place.  In fact, the operation of a power grid can be
viewed as a noisy process.  Rather than attempting to model the noise and
complexity in detail, we seek a generic modeling methodology that can serve to
expose system vulnerabilities.
\end{itemize}
\noindent The approach we take relies on the fact that one can approximate a 
variety of complex physical phenomena that (negatively) affect the performance
of a line by simply perturbing that line's resistance (or, for AC models,
the conductance, susceptance, etc.). In particular, by significantly
increasing the resistance of an arc we will, in general, force the power flow on
that line to zero.  This modeling approach 
becomes particularly effective, from a system perspective, when the resistances
of many arcs are simultaneously altered in an {\em adversarial} fashion. \\

\noindent Accordingly, our second model works as follows:
\begin{itemize}
\item [(I)] The attacker {\em sets} the resistance $x_{ij}$ of any arc $(i,j)$.   
\item [(II)] The attacker is constrained: we must have $x \in F$ for a certain
known set $F$.
\item [(III)] The output of each generator $i$ is fixed at a given value $P_i$,
and similarly each demand value $D_i$ is also fixed at a given value.
\item  [(IV)]  The objective of the attacker is to maximize the overload of
any arc, that is to say, the attacker wants to solve
\begin{eqnarray}
&& \max_{x \in F} \, \max_{ij} \left\{ \frac{ |f_{ij}| }{u_{ij}} \right\}, \label{nonlinobj}
\end{eqnarray}
\noindent where the $f_{ij}$ are the resulting power flows.  
\end{itemize}
\noindent In view of Lemma \ref{unique}, (III) implies that in (d) the vector $f$ is
unique for each choice of $x$; thus the problem is well-posed. \\

\noindent In future work we plan to relax (III).  But (I), (II), (IV) already capture
a great deal of the inherent complexity of power flows. Moreover, suppose that 
e.g. the value of (\ref{nonlinobj}) equals $1.25$. Then even if we allow
demands to be reduced, but insist that this be done under a {\em fair} 
demand-reduction discipline (one that decreases all demands by the same factor)
the system will lose $25 \%$ of the total demand if overloads are to
be avoided (and it is not surprising that the same qualitative 
conclusion holds even if demands are ``unfairly'' reduced to minimize
maximum overload; see Table \ref{comparo}). Thus we expect that the impact of (III), under this model, may not
be severe. \\

\noindent For technical reasons, it will become more convenient to deal with the 
inverses of resistances, the so-called ``conductances.''  For
each $(i,j) \in E$, write $y_{ij} = 1/x_{ij}$, and let $y$ be the vector
of $y_{ij}$.  Likewise, instead of considering a set $F \subseteq \cR^{E}$ 
of allowable 
values $x_{ij}$ we will consider a set $\Gamma$ be a set 
describing the conductance
values that the adversary is allowed to use. We are interested in 
a problem of the form
\begin{eqnarray}
&& \max_{y \in \Gamma} \, \max_{ij} \left\{ \frac{ |f_{ij}(y)| }{u_{ij}} \right\}, \label{nonlinobj2}
\end{eqnarray}
\noindent where as just discussed the notation
$f_{ij}(y)$ is justified.\\

\noindent A relevant example of a set $\Gamma$ is that given by:
\begin{eqnarray}
&& \sum_{ij} \frac{1}{y_{ij}} \, \le \, B,  \ \ \ \ \ \  \frac{1}{x^U_{ij}} \le y_{ij} \le \frac{1}{x^L_{ij}} \ \ \ \ \forall \, (i,j), \label{ybudget}
\end{eqnarray}
\noindent where $B$ is a given 'budget', and, for any arc $(i,j)$, $x^L_{ij}$ and 
$x^U_{ij}$ and 
indicates a minimum and maximum 
value for the resistance at $(i,j)$.  Suppose the initial 
resistances $x_{ij}$ are all equal to some common value $\bar x$, and we set
$x^L_{ij} = \bar x$ for every $(i,j)$, and
$B =  k \, \theta \, \bar x \, + \, (|E| - k ) \bar x$, 
where $k > 0$ is an integer and $\theta > 1$ is large.  Then, roughly speaking,
 we are approximately allowing the adversary to make the
resistance of (up to) $k$ arcs ``very large'',  while not decreasing any resistance,
a problem closely reminiscent of the classical $N - k$ problem.   We will make this 
statement more precise later. 

If the objective in
(\ref{nonlinobj2})  is convex then the optimum will take place at some extreme
point.  In general, the objective is not convex; but 
computational experience shows that we tend to converge to points that are either
extreme points, or very close to extreme points (see the computational section).

\noindent Obviously, the problem we are describing differs from 
the standard N-k problem (though in Section \ref{relationshipnmk} we
present some comparisons). However, in our opinion we obtain a more 
effective approach for handling modeling noise; the capability to handle
much larger problems provides further encouragement.

\subsection{Solution methodology}
Problem (\ref{nonlinobj2}) is not smooth.  However, it is equivalent to:
\begin{eqnarray}
&& \max_{y , p, q} \,\,\,\,\, \sum_{ij} \frac{ f_{ij}(y) }{u_{ij}} ( p_{ij}  - q_{ij} )\label{nonlinobj3} \\
&& \mbox{s.t.} \,\,\,\,\,\,\,\,\,\, \sum_{ij} ( p_{ij} + q_{ij}) \, = \, 1, \label{pqsum}\\
 && \,\,\,\,\,\,\,\,\,\,\,\,\,\,\,\,\,\,\,\, y \in \Gamma, \,\,\,\,\, p, q \ge 0. \label{nonlinsets}
\end{eqnarray} 
\noindent We sketch a
proof of the equivalence.  Suppose $(y^*, p^*, q^*)$ is an optimal solution to (\ref{nonlinobj3})-(\ref{nonlinsets}); let $(\hat i, \hat j)$ be such that
$| f_{\hat i \hat j}(y^*)|/u_{\hat i \hat j} = \max_{ij} | f_{ij}(y^*)|/u_{ij}$. Then without loss of generality if $f_{\hat i \hat j}(y^*) > 0$ 
(resp., $f_{\hat i \hat j}(y^*) \le 0$) we will
have $p^*_{\hat i \hat j} = 1$ (resp., $q^*_{\hat i \hat j} = 1$) and all other $p^*$, $q^*$ equal to zero.  This proves the equivalence once way and the
converse is similar.\\

\noindent In order to work with this formulation we need to develop a more 
explicit representation of the functions $f_{ij}(y)$.  This will require a 
sequence of technical results given in the following section;
however a brief discussion of our approach follows. \\

\noindent Problem (\ref{nonlinobj3}), although smooth, is not concave.  
A relatively recent research thrust has focused on adapting techniques of
(convex) nonlinear programming to nonconvex problems;
resulting in a very large literature with interesting and useful
results; see \cite{fletchetal}, \cite{benshavan}.  Since one is attempting
to solve non-convex minimization (and thus, NP-hard) problems, 
there is no guarantee that a global optimum will be found by these techniques.
One can sometimes assume that a global optimum is approximately known; and the
techniques then are likely to converge to the optimum from an appropriate 
guess.

In any case, (a) the use of nonlinear models allows for much richer
representation of problems, (b) the very successful numerical methodology 
backing convex optimization is brought to bear, and (c) even though only
a local optimum may be found, at least one is relying on an agnostic, 
``honest'' optimization technique as opposed to a pure heuristic or a 
method that makes structural assumptions about the nature of the optimum
in order to simplify the problem.

In our approach we will indeed rely on this methodology -- items (a)-(c) 
precisely capture the reasons for our choice.  Points (a) and (c)
are particularly important in our blackout context: 
we are very keen on modeling the nonlinearities, and on using 
a truly agnostic algorithm to root out hidden weaknesses in a network.
And from a computational perspective, the 
approach does pay off, because we are able to comfortably handle problems
with on the order of 1000 arcs.  

As a final point, note that 
in principle one could rely on a branch-and-bound procedure to 
actually find the global optimum.  This will be a subject for future research.

\subsubsection{Some comments}
As noted above, research such that described in \cite{fletchetal}, \cite{benshavan} has led to effective algorithms that adapt ideas from convex
nonlinear programming to nonconvex settings.  Suppose we consider a 
a linearly constrained problem of the form 
$\min \left\{ \cF(x) \, : \, Ax \ge b \right\}$.  Implementations 
such as LOQO \cite{loqo} or IPOPT \cite{ipopt} require, in addition to some representation of the
linear constraints $Ax \ge b$, subroutines for computing, at any given
point $\hat x$:
\begin{itemize}
\item [(1)] The functional value $\cF(\hat x)$,
\item [(2)] The gradient $\nabla F(\hat x)$, and, ideally,
\item [(3)] The Hessian $\nabla^2 F(\hat x)$. 
\end{itemize}
\noindent If routines for e.g. the computation of the Hessian are not available, 
then automatic differentiation may be used. At each iteration the algorithms will evaluate the subroutines
and perform additional work, i.e. matrix computations. Possibly, 
the cumulative run-time accrued in the computation of (1)-(3) could represent a large fraction of
the overall run-time. Accordingly, there is a premium on developing fast
routines for the three computations given above, especially in large-scale
settings.  This a major goal in the developments given below.

Additionally, a given optimization problem may admit many mathematically
equivalent formulations.  However, different formulations may lead to 
vastly different convergence rates and run-times. This can become 
especially critical in large-scale applications.  Broadly speaking, one 
can seek two (sometimes opposing) goals:
\begin{itemize}
\item [(a)] Compactness, i.e. ``small'' size of a formulation.  This is important
in the sense that numerical linear algebra routines (such as computation 
of Cholesky factorizations) is a very significant ingredient in the algorithms
we are concerned with.  A large reduction in problem size may well lead to
significant reduction in run-times.
\item [(b)] ``Representativity''.  Even if two formulations are equivalent, one
of them may more directly capture the inherent structure of the problem,
in particular, the interaction between the objective function and constraints.
\end{itemize}
\noindent Our techniques achieve {\em both} (a) and (b).  We will construct an 
explicit representation of the functions $f_{ij}(y)$ given above, such that the three
evaluation steps discussed in (1)-(3) indeed admit efficient implementations using sparse
linear algebra techniques.  Moreover, the approach is ``compact'' in that, essentially, the
only variables we deal with are the $y_{ij}$ -- we do not use the straightforward 
indirect representation 
involving not just the $y$ variables, but also variables for the flows and the angles $\theta$.
As we will argue below, our approach does indeed pay off.  Our techniques lead to fast 
convergence, both in terms of the overall run-time and in terms of the iteration count, 
even in cases where the number of lines is on the order of $1000$.\\

\vspace{.1in}

In what follows, we will first provide a review of some relevant material in linear algebra
(Section \ref{laplace}).  This material is used to make some structural remarks 
in Section \ref{observations}. Section \ref{relationshipnmk} presents a result 
relating the model to the standard N-k problem. Section \ref{efficientlinalg} 
describes our algorithms for computing the
gradient and Hessian of the objective function for problem (\ref{nonlinobj3})-(\ref{nonlinsets}).
Finally, Section \ref{impdetails} presents details of our implementation, and Section 
\ref{expers} describes our numerical experiments.

\subsubsection{Laplacians}\label{laplace}

In this section we present some background material on linear algebra and Laplacians of graphs -- the results are 
standard.   See \cite{boyd} for
relevant material. As before we have a connected, directed network $G$ with $n$ nodes and $m$ arcs and with
node-arc incidence matrix $N$.  For
a positive diagonal matrix $Y \in \cR^{m \times m}$ we will write
\begin{eqnarray}
L \, = \, N Y N^{T}, && J \, = \, L \, + \, \frac{1}{n} \, \b1 \b1^{T}. \nonumber
\end{eqnarray}
\noindent where $\b1 \in \cR^n$ is the vector $(1, 1, \ldots, 1)^T$.    $L$ is called a {\em generalized Laplacian}. 
We have
that $L$ is symmetric positive-semidefinite.  If
$\lambda_1 \le \lambda_2 \le \, \ldots \, \le \lambda_n$ are the eigenvalues of
$L$, and $v^1, v^2, \ldots, v^n$ are the corresponding unit-norm eigenvectors, then 
\begin{eqnarray}
\lambda_1 = 0, & \mbox{but} & \lambda_i > 0 \ \ \ \mbox{for} \ i > 1, \nonumber
\end{eqnarray}
\noindent because $G$ is connected, and thus $L$ has rank $n - 1$.   The
same argument shows that since $N \b1 = 0$, we can assume $v^1 = n^{-1/2}~ \b1$.  Finally,
since different eigenvectors are orthogonal, we have $\b1^T v^i \, = \, 0$ 
for $2 \le i \le n$.  Lemmas \ref{landj}, \ref{sol} and \ref{powk} follow from results
in \cite{boyd} and from basic properties of the matrix $N$.

\begin{LE} \label{landj} $L$ and $J$ have the same eigenvectors, and all but one of their eigenvalues coincide. Further, $J$ is invertible.\end{LE}

\begin{LE} \label{sol} Let $b \in \cR^n$. 
Any solution to the system of equations $L \alpha = b$ is of the form
$$ \alpha = J^{-1}b + \delta \b1, $$
for some $\delta \in \cR$. \end{LE} 

\noindent Define $$P = I - J.$$ 
\noindent Note that the eigenvalues of $P$ are $0$ and $1 - \lambda_i$, $2 \le i \le n$;  thus if we have
\begin{eqnarray}
&& \sum_{(u,v)} y_{uv} \, < \, 1/2, \,\,\,\, \mbox{for all $u$},  \label{ybound}
\end{eqnarray}
\noindent then it is not difficult to show that 
\begin{eqnarray}
&& \, 0 \, < \, 1 - \lambda_i < 1, \ \ \ \ \mbox{for all $i \ge 2$}.  \label{lambdabound}
\end{eqnarray}
\noindent (See \cite{mohar} for related background). In such a case we can write
\begin{eqnarray}
J^{-1} \, = \, (I-P)^{-1} \, = \, \sum_{k = 0}^{\infty} P^k. \label{Jinv}
\end{eqnarray}
\begin{LE} \label{powk}
For any integer $k > 0$, $P^k \, = \, (I - N Y N^T)^k - \frac{1}{n}\b1\b1^{T}$. \end{LE} 

\subsubsection{Observations}\label{observations}

Consider problem
(\ref{nonlinobj3})-(\ref{nonlinsets}), where, as per our modeling assumption (III), 
$b$ denote the (fixed) net supply vector, i.e. $b_i = P_i$ for a generator $i$,
$b_i = -D_i$ for a demand node $i$, and $b_i = 0$ otherwise.  Denoting by $Y$ the diagonal matrix with entries $1/y_{ij}$,
we have that given $Y$ the unique power flows $f$ and voltages $\theta$ are obtained by solving the system
\begin{eqnarray}
N^{T}\theta - Y^{-1} f & = & 0 \label{pf1} \\
N f & = & b. \label{pf2}
\end{eqnarray}
\noindent 
In what follows, it will be convenient to assume that condition (\ref{lambdabound}) holds, i.e.
$1 - \lambda_i < 1$ for each $i$.  Next we argue that without loss of generality we
can assume that this holds.

As noted above, this condition will be satisfied if
$\sum_{(u,v)} y_{uv} \, < \, 1/2$ for all $u$ (eq. (\ref{ybound}) above).  Suppose we were
to scale all $y_{i}$ by a common multiplier $\mu > 0$, and instead of system 
(\ref{pf1})-(\ref{pf2}), we consider:
\begin{eqnarray}
N^{T}\theta - \mu^{-1} Y^{-1} f & = & 0 \label{pf1kappa} \\
N f & = & \mu~ b. \label{pf2kappa}
\end{eqnarray}
\noindent We have that $(f, \theta)$ is a solution to (\ref{pf1})-(\ref{pf2}) iff 
$(\mu f, \theta)$ is a solution to (\ref{pf1kappa})-(\ref{pf2kappa}).
Thus, if we assume that the 
set $\Gamma$ in our formulation (\ref{nonlinobj3})-(\ref{nonlinsets}) is bounded
(as is the case if we use (\ref{ybudget})) then, without loss of generality, (\ref{ybound}) 
indeed holds.  Consequently, 
in what follows we will assume that 
\begin{eqnarray}
&& \exists \ r < 1 \ \ \mbox{such that $1 - \lambda_i < r$ for $2 \le i \le n$.} \label{rassumption}   
\end{eqnarray}
\noindent By Lemma \ref{sol} each solution to (\ref{pf1})-(\ref{pf2}) is of the form
\begin{eqnarray}
&& \theta = J^{-1}b + \delta\b1\ \ \ \ \mbox{for some $\delta \in \cR$}, \label{ang}\\
&& f = Y N^{T} J^{-1} b. \nonumber
\end{eqnarray}
\noindent  For each arc $(i,j)$ denote by $c_{ij}$ the column of $N$ corresponding to $(i,j)$. Using (\ref{Jinv}) we therefore have
\begin{eqnarray}
f_{ij} & = & y_{ij} c_{ij}^{T} \left[ I + P + P^2 + P^3+ \ldots \right] b, \ \ \ \forall (i,j), \ \ \ \mbox{and} \label{fexpr} \\
\theta_i - \theta_j \, = \, c_{ij}^{T} \theta & = & c_{ij}^{T}\left [I + P + P^2 + P^3+ \ldots \right] b  \, = \, c_{ij}^{T} \sum_{k = 0}^{\infty} P^{k} ~ b, \nonumber
\end{eqnarray}

\noindent In the following we will be handling expressions with infinite series such as the above. 
In order to facilitate the analysis we need a 'uniform convergence' argument, as follows.  
Given $y \in \Gamma$, note that we can write
\begin{eqnarray}
&& P \, = \, P(y) \, = \, U(y) \Lambda(y) U(y)^T, \nonumber
\end{eqnarray}
\noindent where $U(y)$ is a unitary matrix and $\Lambda(y)$ is the 
diagonal matrix containing the eigenvalues of $P(y)$.  Hence, for any
$k \ge 1$ and any arc $(i,j)$ (and dropping the dependence on $y_{st}$ for simplicity),
\begin{eqnarray}
&& | n^T_{ij} P^k b | \, = \, | n^T_{ij} U \Lambda^k U^T b | \, < \, \nu^k, \label{uniformconv}
\end{eqnarray}
\noindent for some $\nu < 1$, by (\ref{rassumption}).    We will rely on this bound below.

\subsection{Relationship to the standard N-k problem}\label{relationshipnmk}

\noindent As a first consequence of (\ref{uniformconv}) we 
have the following result, showing that appropriate 
assumptions the continuous model we consider
is related to the network vulnerability models in 
Section \ref{nmk}.  

\begin{LE} \label{limitlemma} Let $S$ be a set of arcs 
whose removal does not disconnect $G$.  
Suppose we fix the values $y_{ij} = 1/x_{ij}$ for each arc $(i,j) \notin S$,
and we likewise 
set $y_{st} = \epsilon$ for each arc $(s,t) \in S$. 
Let 
$(f(y), \theta(y))$  denote the resulting power flow on $G$, 
and $(\bar f, \bar \theta)$ the power flow on $G - S$.\\

\noindent Then 
\begin{itemize}
\item[(a)] $ \lim_{\epsilon \rightarrow 0} f_{st}(y) \, = \, 0$, for all $(s,t) \in S$,
\item[(b)] For any $(u,v) \notin S$, $\lim_{\epsilon \rightarrow 0} f_{uv}(y) \, = \, \bar f_{uv}$.
\item[(c)] For any $(u,v)$, $\lim_{\epsilon \rightarrow 0} (\theta_u(y) - \theta_v(y)) \, = \, \bar \theta_u - \bar \theta_v$. 
\end{itemize} \end{LE}
\noindent {\em Proof.} (a) Let $\tilde G = G - (s,t)$, let $\tilde N$ be 
node-arc incidence matrix of $\tilde G$, $\tilde Y$ the restriction of
$Y$ to $E - (s,t)$, and $\tilde P = I - \tilde N \tilde Y \tilde N^T - \frac{1}{n}\b1\b1^{T}$. \\

\noindent For any integer $k \ge 1$ we have by Lemma \ref{powk}
\begin{eqnarray}
&& \lim_{\epsilon \rightarrow0}P^k = \lim_{\epsilon \rightarrow0}(I- N Y N^{T})^k - \frac{1}{n}\b1\b1^{T} = (I - \tilde{N} \tilde{Y} \tilde{N}^{T})^k - \frac{1}{n}\b1\b1^{T} = \tilde{P}^{k}. \nonumber
\end{eqnarray} Consequently, by
(\ref{fexpr}), for any $(s,t) \in S$,
\begin{eqnarray}
\lim_{\epsilon \rightarrow0} f_{st} & = & \lim_{\epsilon \rightarrow0} \left[ y_{st} c_{st}^T \left( \sum_{k=0}^{\infty} P^k \right) b\right] \,= \,  \sum_{k=0}^{\infty}\left[  \lim_{\epsilon\rightarrow0}  y_{st}\left( c_{st}^T P^k b \right) \right] \, = \, 0, \nonumber
\end{eqnarray}
\noindent where the exchange between summation and limit is valid because of
(\ref{uniformconv}).   The proof of (b), (c) are similar.\QED \\

\noindent Lemma \ref{limitlemma} can be interpreted as describing a particular 
attack pattern -- make $x_{ij}$ very large for $(i,j) \in S$ and leave all 
other $x_{ij}$ unchanged. Our computational experiments
show that the pattern assumed by the Lemma
is approximately correct: given an attack budget, the attacker tends to concentrate most
of the attack on a small number of arcs (essentially, making their resistance very large),
while at the same time attacking a larger number of lines with a small portion of the budget.

\subsection{Efficient computation of the gradient and Hessian}\label{efficientlinalg}

\noindent In the following set of results we determine efficient closed-form expressions for
the gradient and Hessian of the objective in (\ref{nonlinobj2}).  As before, we denote 
by $c_{ij}$ the column of the node-arc incidence matrix of the network corresponding to
arc $(i,j)$.  First we present a technical result.  This will be followed 
by the development of formulas for the gradient (eqs. (\ref{grad-1})-(\ref{grad-2})) and the Hessian (eqs. (\ref{hess-1})-(\ref{hess-3})).

\begin{LE} \label{partial1} For any integer $k > 0$, and any arc $(i,j)$
\begin{eqnarray}
&(a)& \b1^{T} P^{k} \, = \, 0, \nonumber \\
&(b)& \frac{\partial}{\partial{y_{ij}}} \left[P^k b \right] \, = \, P \frac{\partial}{\partial{y_{ij}}}\left[P^{k-1} b \right] - c_{ij} c_{ij}^{T}P^{k-1}b. \nonumber
\end{eqnarray}
\end{LE}

\noindent {\em Proof.} Note that $\b1^{T}P = \b1^{T}(I - J) = \b1^{T}(I- NYN^{T}- \frac{1}{n}\b1\b1^{T}) = 0$. Hence $\b1^{T}P^{k} = 0 $.  

\begin{eqnarray*}
\frac{\partial}{\partial{y_{ij}}}\left[P^k b \right] &=& 
\frac{\partial}{\partial{y_{ij}}}\left[PP^{k-1} b \right] \\
&=& \frac{\partial}{\partial{y_{ij}}}\left[\left(I - \sum_{ (u,v) \in E} y_{uv} \, c_{uv} c_{uv}^{T}- \frac{1}{n} \b1\b1^{T} \right) P^{k-1} b \right]\\
&=& \frac{\partial}{\partial{y_{ij}}}\left[P^{k-1} b \right] - \frac{\partial}{\partial{y_{ij}}}\left[\left(\sum_{(u,v) \in E} y_{uv} \, c_{uv} c_{uv}^{T} \right)P^{k-1} b \right]- \frac{\partial}{\partial{y_{ij}}}\left[\frac{1}{n} \b1\b1^{T}P^{k-1} b  \right] \\
&=& \frac{\partial}{\partial{y_{ij}}}\left[P^{k-1} b \right] - \frac{\partial}{\partial{y_{ij}}}\left[\left(\sum_{(u,v) \in E} y_{uv} \, c_{uv} c_{uv}^{T} \right)P^{k-1} b \right]\\
&=& \frac{\partial}{\partial{y_{ij}}}\left[P^{k-1} b \right] - \sum_{(u,v) \in E}  \frac{\partial}{\partial{y_{ij}}}\left[y_{uv} \, c_{uv} c_{uv}^{T}P^{k-1} b \right]\\
&=& \frac{\partial}{\partial{y_{ij}}}\left[P^{k-1} b \right] - \sum_{(u,v) \in E}  \left[\frac{\partial y_{uv} \,}{\partial{y_{ij}}} \right] \ \ c_{uv} c_{uv}^{T}P^{k-1} b
- \sum_{(u,v) \in E}  y_{uv} \, \frac{\partial}{\partial{y_{ij}}} \left[c_{uv} c_{uv}^{T}P^{k-1} b \right]\\
&=& \frac{\partial}{\partial{y_{ij}}}\left[P^{k-1} b \right] - c_{ij} c_{ij}^{T}P^{k-1} b
- \sum_{(u,v) \in E}  y_{uv} \, c_{uv} c_{uv}^{T} \frac{\partial}{\partial{y_{ij}}} \left[P^{k-1} b \right]\\
&=& \left[I - \sum_{(u,v) \in E}  y_{uv} \, c_{uv} c_{uv}^{T} \right] \frac{\partial}{\partial{y_{ij}}}\left[P^{k-1} b \right] - c_{ij} c_{ij}^{T}P^{k-1} b\\
&=& \left[P + \frac{1}{n} \b1\b1^{T} \right] \frac{\partial}{\partial{y_{ij}}}\left[P^{k-1} b \right] - c_{ij} c_{ij}^{T}P^{k-1} b\\
&=& P \frac{\partial}{\partial{y_{ij}}}\left[P^{k-1} b \right] - c_{ij} c_{ij}^{T}P^{k-1}b +  \frac{\partial}{\partial{y_{ij}}}\left[\frac{1}{n} \b1\b1^{T}P^{k-1} b \right]\\
&=& P \frac{\partial}{\partial{y_{ij}}}\left[P^{k-1} b \right] - c_{ij} c_{ij}^{T}P^{k-1}b.\\
\end{eqnarray*}
\noindent where the third and the last equality follow from (a). \QED\\

\noindent Using the above recursive formula we can write the following expressions:
\begin{eqnarray*}
\frac{\partial}{\partial{y_{ij}}}[P b] &=& -c_{ij} c_{ij}^{T}b \nonumber\\
\frac{\partial}{\partial{y_{ij}}}[P^{2}b] &=& P\frac{\partial}{\partial{y_{ij}}}[P b] - c_{ij} c_{ij}^{T}P b \\
\frac{\partial}{\partial{y_{ij}}}[P^{3}b] &=& P^{2}\frac{\partial}{\partial{y_{ij}}}[P b] - P c_{ij} c_{ij}^{T}P b - c_{ij} c_{ij}^{T}P^{2}b\\
\frac{\partial}{\partial{y_{ij}}}[P^{4}b] &=& P^{3}\frac{\partial}{\partial{y_{ij}}}[P b] - P^{2} c_{ij} c_{ij}^{T}P b - P c_{ij} c_{ij}^{T}P^{2}b - c_{ij} c_{ij}^{T}P^{3}b\\
&\vdots&\\
\frac{\partial}{\partial{y_{ij}}}[P^{k}b] &=& P^{k-1}\frac{\partial}{\partial{y_{ij}}}[P b] - P^{k-2} c_{ij} c_{ij}^{T}P b - P^{k-3}c_{ij} c_{ij}^{T}P^{2}b - \ldots - c_{ij} c_{ij}^{T}P^{k-1}b\\
\end{eqnarray*}

\noindent Consequently,  defining
\begin{eqnarray}
\widetilde{\nabla}_{ij} & = & \frac{\partial}{\partial{y_{ij}}} \left[I + P + P^2 + \ldots \right] b, 
\end{eqnarray}
\noindent we have
\begin{eqnarray}
\widetilde{\nabla}_{ij}  &=& \left[I + P + P^2 + \ldots \right]\frac{\partial}{\partial{y_{ij}}}\left[P b \right] \, - \, \left(I + P + P^2 + \ldots \right) \, c_{ij} c_{ij}^{T} \, \left(P + P^2 + P^3+ \ldots \right)b \nonumber\\
&=& -\left[I + P + P^2 + \ldots \right]c_{ij} c_{ij}^{T}b \, - \, \left(I + P + P^2 + \ldots \right)c_{ij} c_{ij}^{T}\left(I + P + P^2 + \ldots - I  \right)b \nonumber\\
&=& - \left(I + P + P^2 + \ldots \right)c_{ij} c_{ij}^{T}\left(I + P + P^2 + \ldots  \right)b \nonumber\\
&=& -J^{-1} \, c_{ij} c_{ij}^{T} \, \theta, \nonumber
\end{eqnarray}
\noindent where the last equality follows from (\ref{ang}) and (\ref{Jinv}),
and the fact that $c_{ij}^{T} \b1 = 0$.\\

\noindent Using (\ref{fexpr}), the gradient of function $f_{uv}(y)$ with respect to the variables $y_{ij}$ can be written as:
\begin{eqnarray}
 \frac{\partial{f_{uv}}}{\partial{y_{ij}}} &=& y_{uv} \, c_{uv}^{T} \, \frac{\partial}{\partial{y_{ij}}}\left[I + P + P^2 + P^3+ \ldots \right]b \, = \, 
y_{uv} \, c_{uv}^{T} \, \widetilde{\nabla}_{ij},  \ \ \ \ (i,j) \neq (u,v) \label{grad-1}\\
\frac{\partial{f_{ij}}}{\partial{y_{ij}}} &=& c_{ij}^{T}\left[I + P + P^2 + P^3+ \ldots \right]b
\, + \, y_{ij} \, c_{ij}^{T} \, \frac{\partial}{\partial{y_{ij}}}\left[I + P + P^2 + P^3+ \ldots \right]b \nonumber \\
&=& c_{ij}^{T} \widetilde{\nabla}_{ij} \, + \, y_{ij} c_{ij}^{T} \widetilde{\nabla}_{ij}. \label{grad-2}
\end{eqnarray}

\noindent We similarly develop closed-form expressions for the second-order derivatives. 
For $(u,v) \neq (i,j),  (u,v) \neq (h,k)$, we have the following :
\begin{eqnarray}
 \frac{\partial^{2}{f_{uv}}}{\partial{y_{ij}} \partial{y_{hk}}} &=& y_{uv} c_{uv}^{T} \, [ \, (I + P + P^2 + P^3 + \ldots) \, c_{ij} c_{ij}^{T} \, (I + P + P^2 + P^3 + \ldots)
 c_{hk} c_{hk}^{T} \nonumber \\
 &&  + \, (I + P + P^2 + P^3+ \ldots) \, c_{hk} c_{hk}^{T} \, (I + P + P^2 + P^3+ \ldots) \, c_{ij} c_{ij}^{T}  \,] \theta  \nonumber \\
 &=& - y_{uv} c_{uv}^{T} J^{-1} \, \left[ c_{ij} c_{ij}^{T} \widetilde{\nabla}_{hk} + c_{hk} c_{hk}^{T} \widetilde{\nabla}_{ij}  \right]. \label{hess-1}
\end{eqnarray}

\noindent Similarly, the remaining terms are:
\begin{eqnarray}
 \frac{\partial^{2}{f_{uv}}}{\partial{y_{uv}^{2}}} &=&  2 \, c_{uv}^{T} \widetilde{\nabla}_{uv} \,\, - \,\,2 \, y_{uv} c_{uv}^{T} J^{-1} c_{uv} c_{uv}^{T} \widetilde{\nabla}_{uv}, \label{hess-2}\\
 \frac{\partial^{2}{f_{uv}}}{\partial{y_{uv}} \partial{y_{ij}}} &=& c_{uv}^{T} \, \widetilde{\nabla}_{ij} \,\, - \,\, y_{uv} \,c_{uv}^{T} J^{-1} \left[ c_{ij} c_{ij}^{T} \, \widetilde{\nabla}_i \, + \, c_{uv} c_{uv}^{T} \, \widetilde{\nabla}_{ij}  \right] \label{hess-3}
\end{eqnarray}

\subsection{Implementation details}\label{impdetails}
We use LOQO \cite{loqo} to solve problem (\ref{nonlinobj3})-(\ref{nonlinsets}), using
$\Gamma \, = \, \left\{ \, y \ge 0 \, : \sum_{ij} \frac{1}{y_{ij}} \, \le \, B \, \right\}$
with values of $B$ that we selected. LOQO is an infeasible primal-dual, interior-point method 
applied to a sequence of quadratic approximations
to the given problem. The procedure stops if at any iteration the 
primal and dual problems are feasible and with objective values that are $\it{close} $ 
to each other, in which case a local optimal solution is found. 
Additionally, LOQO uses an upper bound on the overall
number of iterations it is allowed to perform. \\

\noindent As outlined above, each iteration of LOQO 
requires the Hessian and gradient of the objective function. We do this as 
described above, e.g. using (\ref{grad-1}), (\ref{grad-2}), (\ref{hess-1})-(\ref{hess-3}).
This approach requires the computation of quantities
$c_{uv}^T J^{-1} c_{ij}$ for each pair of arcs $(i,j)$, $(u,v)$.
At any given iteration, we compute and store these 
quantities  (which can be done in  $O(n^{2} + n m)$ space). 

\noindent In order to compute $c_{uv}^T J^{-1} c_{ij}$, for given $(i,j)$ and $(u,v)$, 
we simply solve the sparse
linear system on variables $\kappa$, $\lambda$:
\begin{eqnarray}
N^{T}\kappa - Y^{-1} \lambda & = & 0 \label{pfaux} \\
N \lambda & = & c_{ij}. \label{pfaux2}
\end{eqnarray}
\noindent As in (\ref{ang}), we have $\kappa = J^{-1}c_{ij} + \delta\b1$ for
some real $\delta$.  But then $c_{uv}^T \kappa = c_{uv}^T J^{-1}c_{ij}$, 
the desired quantity.  In
order to solve (\ref{pfaux})-(\ref{pfaux2}) we use Cplex 
(though any efficient linear algebra package would suffice).\\

\noindent We point out that, alternatively, LOQO (as well as other solvers
such as IPOPT) can perform symbolic
differentiation in order to directly compute the Hessian and gradient.  We 
could in principle follow this approach in order to solve a problem
with objective (\ref{nonlinobj3}), constraints (\ref{pqsum}), (\ref{nonlinsets}) {\em and} (\ref{1b}), (\ref{ohm-eq}). We prefer our approach 
because it employs fewer variables (we do not need the flow variables 
or the angles) and
primal feasibility is far simpler.  This is the ``compactness'' ingredient alluded to before.\\

\noindent In our implementation,
we fix a value for the iteration limit, but apply additional stopping criteria:
\begin{itemize}
\item[(1)] If both primal and dual are feasible, 
we consider the relative error between the primal 
and dual values, $\epsilon = \frac{\mbox{PV - DV}}{\mbox{DV}}$, where 'PV' 
and 'DV' refer to primal and dual values respectively. If the relative error 
$\epsilon$ is less than some desired threshold we stop, and report the
solution as ``$\epsilon$-locally-optimal.'' 
\item[(2)] If on the other hand we reach the iteration limit without 
stopping, then we consider the
last iteration at which we had both primal and dual feasible solutions. 
If such an iteration exists, then we report the corresponding configuration
of resistances along with the associated maximum congestion value. If such an iteration 
does not exist, then the report the run as unsuccessful.
\end{itemize}

\noindent Finally, we provide to LOQO the starting point $x_{ij} = x_{ij}^L$
for each arc $(i,j)$.

\subsection{Experiments}\label{expers}
In the experiments reported in this section we used a 2.66 GHz Xeon
machine with 2 MB L2 cache and 16 GB RAM.  All experiments were run using
a single core. Altogether we report on $37$ runs
of the algorithm.  

\subsubsection{Data sets}

For our tests we used the 58- and 118-bus test cases as in 
Section \ref{comp-model1} with some variations on the capacities; as well as the $49$-node ``square grid'' example
and three larger networks
 created using the replication technique described at the start of
Section \ref{comp-model1}: a $300$-node, $409$-arc network, a 
$600$-node, $990$-arc network, and a 
$619$-node, $1368$-arc network.  Additional artificial networks were
created to test specific conditions.  All data sets are available for
download (http://www.columbia.edu/$\sim $dano/research/pgrid/Data.zip). 

We considered several three constraint sets $\Gamma$ as in (\ref{ybudget}):
\begin{itemize}
\item [(1)] \bmath{\Gamma(1)}, where for all $(i,j)$, $x_{ij}^L = 1$ and $x_{ij}^U = 5$,
\item [(2)] \bmath{\Gamma(2)}, where for all $(i,j)$, $x_{ij}^L = 1$ and $x_{ij}^U = 10$,
\item [(3)] \bmath{\Gamma(3)}, where for all $(i,j)$, $x_{ij}^L = 1$ and $x_{ij}^U = 20$.
\end{itemize}
\noindent In each case, we set $B = \sum_{(i,j)}x_{ij}^{L} \, + \, \Delta \mbox{B}$,
where $\Delta \mbox{B}$ represents an ``excess budget''. Note that for example in the case $\Delta \mbox{B} = 30$, under $\Gamma(2)$,  
the attacker can increase (from their minimum value)
the resistance of up to $3$ arcs by a factor of $10$ (with $3$ units of
budget left over).  And under $\Gamma(1)$, up to $6$ arcs can have
their resistance increased by a factor of $5$.  In either case we
have a situation reminiscent of the $N - k$ problem, with small $k$. 

\subsubsection{Focus of the experiments}
In these experiments, we first study how the algorithm 
scales as network size increases (up to on other order of $1000$) and
as $\Delta B$ increases.  A second point of focus is the stability of
the underlying (nonlinear) solver -- e.g., does our algorithm frequently
produce poor results because the solver experiences numerical difficulties.

Next, is there significant impact of alternate
starting point choices for the algorithm, and does that constitute 
evidence of lack of robustness.

An important point we want to study concerns the {\em structure} of the
solutions produced by the algorithm -- what is the 
distribution of the $x_{ij}$ obtained at termination, and is there a logic
to that distribution? 

A final set of experiments carry out a comparison
with results obtained using the standard $N-k$ model.

\subsubsection{Basic run behavior}
Tables \ref{57bus2}-\ref{gamma2large} present results for
different networks and scenarios.  Each column corresponds to a different
value of \bmath{\Delta B}.  For each run, {\bf ``Max Cong''} is the numerical value of
the maximum arc congestion 
(as in (\ref{nonlinobj})) at termination. Additionally, we present the
CPU time (in seconds) taken by the algorithm, the number of iterations, and
the termination criterion, which is indicated by ``Exit Status'', with
the following interpretation:
\begin{itemize}
\item[(1)] '\bmath{\epsilon}{\bf-L-opt.}': the algorithm computed an $\epsilon$-locally-optimal solution.
\item[(2)] '{\bf PDfeas, Iter: lastItn}': the algorithm 
reached the iteration limit without finding an $\epsilon$-locally-optimal solution, but there was an iteration at which both primal 
and dual problems 
were feasible. 'lastItn' gives the last iteration at 
which both primal and dual 
solutions were feasible. 
\item[(3)] '{\bf opt.}': the algorithm attained LOQO's internal optimality
tolerance.
\end{itemize}

\noindent Tables \ref{57bus2} and \ref{118bus2} contain results 
for the 57- and 118-bus networks, respectively, both using set $\Gamma(2)$. 
 Tables \ref{gamma1small} and \ref{gamma2small} handle 
the $49$-node, $84$-arc network, with $14$ demand nodes and 
$4$ generators that we considered
in section \ref{comp-model1}, using sets $\Gamma(1)$ and $\Gamma(2)$ 
respectively. 

\begin{table}[ht]
\caption{\textbf{\emph{57 nodes, 78 arcs, constraint set \bmath{\Gamma(2)} }}}
\label{57bus2}
\hskip 12 pt
\centering
\mbox{\bf{Iteration Limit}}: \bmath{700 , \, \, \epsilon = 0.01}\\
\begin{tabular}{|c|c|c|c|c|}
\hline
\multicolumn{1}{|c|}{}  & \multicolumn{4}{|c|}{ \bmath{\Delta}\bf{B}} \\
\multicolumn{1}{|c|}{}  &  \multicolumn{1}{c}{\bf{9}} &  \multicolumn{1}{c}{\bf{18}} & \multicolumn{1}{c}{\bf{27}} & \multicolumn{1}{c|}{\bf{36}} \\
\hline \hline
   & &  &  &    \\
{\bf Max Cong} & 1.070  & 1.190 & 1.220 & 1.209 \\
   & &  &  &    \\
{\bf Time (sec)} & 8 & 19 & 19 & 19   \\
   & &  &  &    \\
{\bf Iterations} & 339 & Limit & Limit & Limit   \\
   & &  &  &    \\
{\bf Exit Status} & $\epsilon$-L-opt. & PDfeas. & PDfeas. & PDfeas.   \\
   &  & Iter: 700 & Iter: 700 &  Iter: 700  \\
\hline
\end{tabular}
\end{table}

\begin{table}[ht]
\caption{\textbf{\emph{118 nodes, 186 arcs, constraint set \bmath{\Gamma(2)} }}}
\label{118bus2}
\hskip 12 pt
\centering
\mbox{\bf{Iteration Limit}}: \bmath{700 , \, \, \epsilon = 0.01}\\
\begin{tabular}{|c|c|c|c|c|}
\hline
\multicolumn{1}{|c|}{}  & \multicolumn{4}{|c|}{ \bmath{\Delta}\bf{B}} \\
\multicolumn{1}{|c|}{}  &  \multicolumn{1}{c}{\bf{9}} &  \multicolumn{1}{c}{\bf{18}} & \multicolumn{1}{c}{\bf{27}} & \multicolumn{1}{c|}{\bf{36}} \\
\hline \hline
   & &  &  &    \\
{\bf Max Cong} & 1.807  & 2.129 & 2.274 & 2.494 \\
   & &  &  &    \\
{\bf Time (sec)} & 88 & 200 & 195 & 207   \\
   & &  &  &    \\
{\bf Iterations} & Limit & 578 & Limit &  Limit \\
   & &  &  &    \\
{\bf Exit Status} & PDfeas. & $\epsilon$-L-opt. & PDfeas. & PDfeas.   \\
   & Iter: 302 &  & Iter: 700 & Iter: 700    \\
\hline
\end{tabular}
\end{table}

\begin{table}[ht]
\caption{\textbf{\emph{49 nodes, 84 arcs, constraint set \bmath{\Gamma(1)} }}} 
\label{gamma1small}
\hskip 12 pt
\centering
\mbox{\bf{Iteration Limit}}: \bmath{800 , \, \, \epsilon = 0.01}
\begin{tabular}{|c|c|c|c|c|c|c|}
\hline
\multicolumn{1}{|c|}{}  & \multicolumn{6}{|c|}{ \bmath{\Delta}\bf{B}} \\
\multicolumn{1}{|c|}{}  & \multicolumn{1}{c}{\bf{5}} & \multicolumn{1}{c}{\bf{10}} & \multicolumn{1}{c}{\bf{15}} & \multicolumn{1}{c}{\bf{20}} & \multicolumn{1}{c}{\bf{25}} & \multicolumn{1}{c|}{\bf{30}} \\
\hline \hline
   & &  &  &   &  & \\
{\bf Max Cong} & 0.673054 & 0.750547 & 0.815623 & 0.865806 & 0.901453 & 0.951803 \\
   & &  &  &   &  & \\
{\bf Time (sec)} & 12 & 15 & 18 & 19 & 28 & 22 \\
   & &  &  &   &  & \\
{\bf Iterations} & 258 & 347 & 430 & 461 & Limit & 492 \\
   & &  &  &   &  & \\
{\bf Exit Status} & $\epsilon$-L-opt. & $\epsilon$-L-opt. & $\epsilon$-L-opt. & $\epsilon$-L-opt. & PDfeas & $\epsilon$-L-opt. \\
 & & &  & &Iter: 613 & \\
\hline
\end{tabular}

\caption{\textbf{\emph{49 nodes, 84 arcs, constraint set \bmath{\Gamma(2)} }}} 
\label{gamma2small}
\hskip 12 pt
\centering
\mbox{\bf{Iteration Limit}}: \bmath{800 , \, \, \epsilon = 0.01}
\begin{tabular}{|c|c|c|c|c|c|c|}
\hline
\multicolumn{1}{|c|}{}  & \multicolumn{6}{|c|}{ \bmath{\Delta}\bf{B}} \\
\multicolumn{1}{|c|}{}  & \multicolumn{1}{c}{\bf{5}} & \multicolumn{1}{c}{\bf{10}} & \multicolumn{1}{c}{\bf{15}} & \multicolumn{1}{c}{\bf{20}} & \multicolumn{1}{c}{\bf{25}} & \multicolumn{1}{c|}{\bf{30}} \\
 \hline \hline
   & &  &  &   &  & \\
 {\bf Max Cong} & 0.67306 & 0.751673 & 0.815584 & 0.8685 & 0.91523 & 0.9496   \\
   & &  &  &   &  & \\
 {\bf Time (sec)} & 9 & 13 & 34 & 3 & 29 & 30   \\
   & &  &  &   &  & \\
 {\bf Iterations} & 177 & 295 & Limit & Limit & Limit & Limit   \\
   & &  &  &   &  & \\
 {\bf Exit Status} & $\epsilon$-L-opt. & $\epsilon$-L-opt. & PDfeas & PDfeas & PDfeas & PDfeas   \\
 &  &  & Iter: 800 & Iter: 738 & Iter: 624 & Iter: 656   \\
\hline
\end{tabular}
\end{table}

Table \ref{gamma2med} presents similar results for the network with 
$300$ nodes, $409$ arcs ($42$ generators and $172$ loads). 
Note that for the runs $\Delta \mbox{B} \geq 20$ 
the maximum load value is identical; the optimal solution values $x_{ij}$ 
were nearly identical, independent of the initial point given to LOQO.
\begin{table}[ht]
\caption{\textbf{\emph{300 nodes, 409 arcs, constraint set \bmath{\Gamma(2)} }}}
\label{gamma2med}
\hskip 12 pt
\centering
\begin{tabular}{|c|c|c|c|c|}
\multicolumn{5}{c}{\bf{Iteration Limit}: \bmath{500 , \, \, \epsilon = 0.01}}\\
\hline
\multicolumn{1}{|c|}{}  & \multicolumn{4}{|c|}{ \bmath{\Delta}\bf{B}} \\
\multicolumn{1}{|c|}{}  &  \multicolumn{1}{c}{\bf{9}} &  \multicolumn{1}{c}{\bf{18}} & \multicolumn{1}{c}{\bf{27}} & \multicolumn{1}{c|}{\bf{36}} \\
\hline \hline
   & &  &  &    \\
{\bf Max Cong} & 0.590690 & 0.694101 & 0.771165 & 0.771165   \\
   & &  &  &    \\
{\bf Time (sec)} & 208 & 1248 & 981 & 825   \\
   & &  &  &    \\
{\bf Iterations} & 91 & Limit & 406 & 320  \\
   & &  &  &    \\
 {\bf Exit Status} & opt. & PDfeas & opt. & opt  \\
  &  & Iter: 318 &  &      \\
\hline
\end{tabular}
\end{table}

\begin{table}[h]
\caption{\textbf{\emph{600 nodes, 990 arcs, constraint set \bmath{\Gamma(2)} }}}
\label{gamma2large}
\hskip 12 pt
\centering
\mbox{\bf{Iteration Limit}}: \bmath{300 , \, \, \epsilon = 0.01}
\begin{tabular}{|c|c|c|c|c|c|}
\hline
\multicolumn{1}{|c|}{}  & \multicolumn{5}{|c|}{ \bmath{\Delta}\bf{B}} \\
\multicolumn{1}{|c|}{}  &  \multicolumn{1}{c}{\bf{10}} &  \multicolumn{1}{c}{\bf{20}} & \multicolumn{1}{c}{\bf{27}} & \multicolumn{1}{c}{\bf{36}} & \multicolumn{1}{c|}{\bf{40}}\\
\hline \hline
   & &  &  &  &  \\
{\bf Max Cong} & 0.082735 (0.571562) & 1.076251 & 1.156187 & 1.088491  & 1.161887\\
   & &  &  &   & \\
{\bf Time (sec)} & 11848 & 7500 & 4502 & 11251  & 7800\\
   & &  &  &   & \\
{\bf Iterations} & Limit & 210 & 114 & Limit  & 208\\
   & &  &  &   & \\
{\bf Exit Status} & PDfeas & $\epsilon$-L-opt. & $\epsilon$-L-opt. &  PDfeas & $\epsilon$-L-opt. \\
   & Iter: 300 &  &  &  Iter: 300  &\\
\hline
\end{tabular}
\end{table}

Table \ref{gamma2large} contains the results for the network with $600$ nodes, 
and $990$ arcs ($344$ 
demand nodes and $98$ generators) under set $\Gamma(2)$. 
We observed an interesting issue in the case where $\Delta \mbox{B} = 10$. Here, 
LOQO terminated with a solution in which for some arc $(i,j)$, both $p_{ij} > 0$ and $q_{ij} > 0$ (refer to formulation (\ref{nonlinobj3})-(\ref{nonlinsets})). 
The value in parenthesis indicates the true value of the maximum congestion obtained by solving the network 
controller's problem if we were to use the resistance values  ($x_{ij}$) given by LOQO.

\begin{table}[h]
\caption{\textbf{\emph{649 nodes, 1368 arcs, constraint set \bmath{\Gamma(2)}}}}
\label{gamma2huge}
\hskip 12 pt
\centering
\mbox{{\bf Iteration Limit: 500, \bmath{\epsilon = 0.01}}}
\begin{tabular}{|c|c|c|c|c|}
\hline
\multicolumn{1}{|c|}{}  & \multicolumn{4}{|c|}{ \bmath{\Delta}\bf{B}} \\
\multicolumn{1}{|c|}{}  &   \multicolumn{1}{c}{\bf{20}} & \multicolumn{1}{c}{\bf{30}} & \multicolumn{1}{c}{\bf{40}} & \multicolumn{1}{c|}{\bf{60}}\\
\hline \hline
   &  &  &   & \\
{\bf Max Cong} & (0.06732) 1.294629 & 1.942652 & (0.049348) 1.395284 & 2.045111\\
   &  &  &   & \\
{\bf Time (sec)} &  66420 & 36274 & 54070  & 40262\\
   &  &  &   & \\
{\bf Iterations} & Limit & 374 & Limit  & Limit\\
   &  &  &   & \\
{\bf Exit Status} &  DF & $\epsilon$-L-opt. &  DF & PDfeas \\
   &  &  &   & Iter: 491\\
\hline
\end{tabular}
\end{table}

Finally, Table \ref{gamma2huge} presents experiments on the network with
$649$ nodes and $1368$ arcs.  Here, exit status 'DF' means that dual feasibility
was achieved, but not primal feasibility.  In such a case, the
budget constraint 
(\ref{ybudget}) was violated -- the largest (scaled) violation we observed was
$1e-03$.  Even though this is a small violation, 
LOQO's threshold for primal feasibility is $1e-06$; we simply scaled down any resistance value
$x_{ij} > x_{ij}^{min}$ so as to obtain a solution satisfying (\ref{ybudget}).
In Table \ref{gamma2huge}, the quantity following the parenthesis in the 
``Max Cong'' line indicates the resulting maximum congestion, obtained by
solving a controller's problem on the network using the reduced resistance values.\\

{\bf Comments:} The algorithm appears to scale, fairly reliably, to 
cases with approximately $1000$ arcs; at that point the internal solver 
(LOQO) starts to develop some difficulties.  

For any given network, note that the computed solution does vary as a
function of the parameter $\Delta B$, and in the expected manner,
as reflected by the ``Max Cong'' values. However the performance
of the algorithm (running time or number of iterations) 
appears stable as a function of $\Delta B$.  By ``stable'' what we mean
is that
even though larger $\Delta B$
values correspond to larger numbers of arcs that could be maximally
interdicted, the workload incurred by the algorithm {\em does not} increase
``combinatorially'' function of $\Delta B$.  In our opinion, this is
a significant distinction between this algorithm and the algorithm 
presented above for the $N - k$ problem. 

\begin{figure}[htb] 
\centering
\includegraphics[width=4in, angle=270]{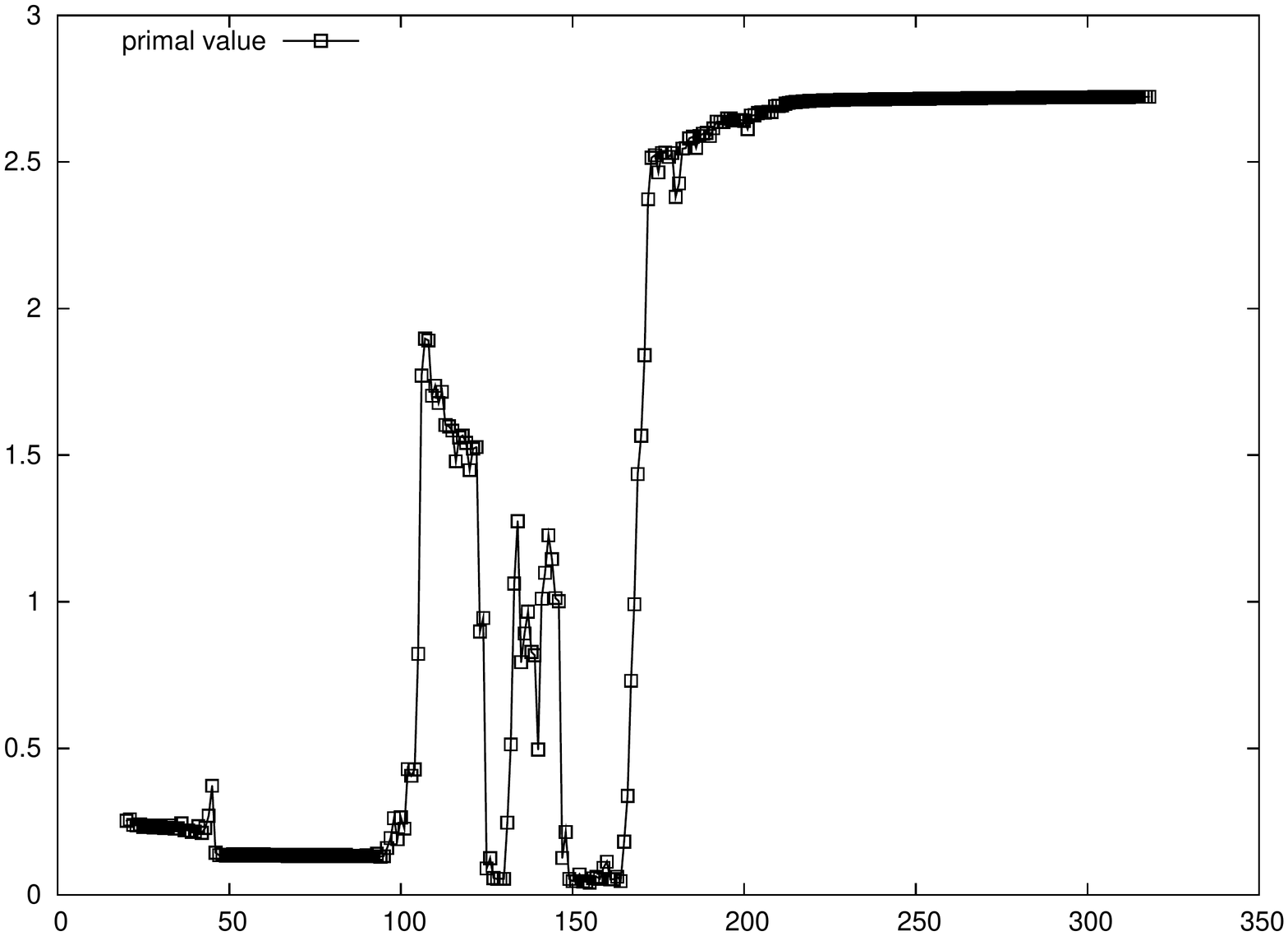}
\caption{Primal values approaching termination.\label{fig:primal}}
\end{figure}

To put it differently, the algorithm in this section appears to allow for
practicable analyses of the impact of multiple choices $\Delta B$; this 
is a critical feature in that  parameterizes the risk-aversion
of the model.  \\

Figure \ref{fig:primal} presents a different view on the progress on
a typical run.  This run concerns the network in Table
\ref{gamma2large} (600 nodes, 990 arcs).  The chart shows the primal
value computed by LOQO on the last 299 iterations (shown on the y-axis; the
x-axis displays the iterations).  It appears that the algorithm computes
several local optima and then settles for a long hill climb.

\subsubsection{Alternative starting points}

The question we consider here is how the final solution computed
by the algorithm varies as a function of the starting point.

Table \ref{gamma2alternative} shows runs using
the network with $49$ nodes and $90$ arcs, using set $\Gamma(3)$ with
$\Delta B = 57$.  For each run we list the maximum congestion at 
termination, and the top six arcs interdicted arcs, with the corresponding
resistance values in parentheses.  Four different choices of starting 
point were considered.

For the first test the starting point was constructed by setting all resistances values to their
lower bounds, i.e., $x_{ij} = 1$. For the second,
we set the resistance of three randomly selected arcs to the
maximum value,
while the remaining arcs were set to the lower bound, 
i.e $x_{ij} = 20, (i,j) \in I \subset E, |I| = 3, \  x_{kl} = 1, {(k,l)} \in E \setminus I $.

For the third test, we set the resistance of six randomly chosen arcs 
to half of the maximum, 
and the resistance of the remaining arcs were set to the minimum value, i.e, 
$x_{ij} = 10, (i,j) \in I \subset E, |I| = 6, \ x_{kl} = 1, (k,l) \in E \setminus I $. For the last test, we used, as starting point, 
the solution found in the test using the 
third starting point.

We note that there is a difference of (at most) $1.5 \%$ in the 
resulting congestion value; while at the same time, and more crucially, the
{\em set} of heavily interdicted arcs does not change. Results such as these
are typical of what we have found in our experiments.  

\begin{table}[htb]
\caption{\textbf{\emph{Impact of changing the starting point}}}
\label{gamma2alternative}
\hskip 12 pt
\centering
\begin{tabular}{|c|c|c|c|c|}
\hline
\multicolumn{1}{|c|}{}  & \multicolumn{4}{|c|}{ {\bf Test} } \\
\multicolumn{1}{|c|}{}  &  \multicolumn{1}{c}{\bf{1}} &  \multicolumn{1}{c}{\bf{2}} & \multicolumn{1}{c}{\bf{3}} & \multicolumn{1}{c|}{\bf{4}} \\
 \hline \hline
  & &  &  &    \\
{\bf Max Cong} & 2.149673 & 2.127635 & 2.164906 & 2.181400   \\
   & &  &  &    \\
   & 29(7.79), 27(7.20) & 29(7.79), 27(7.23)& 29(8.73), 27(8.21) & 29(8.37), 27(7.80) \\
{\bf Top $6$ Arcs} & 41(7.03), 67(7.02) & 41(6.91), 67(7.97) & 41(7.03), 67(7.02) & 41(7.57), 67(7.54)   \\
   & 54(6.72), 79(5.71) & 54(6.58), 79(5.53) & 54(7.52), 79(6.48) & 54(7.24), 79(6.26)   \\
\hline
\end{tabular}
\end{table}
\hspace{.1in} \newline

A conclusive determination of whether our algorithm finds solutions which in some sense are ``close''
to a global optimum is a project for future research which we are now 
undertaking.

\subsubsection{Distribution of attack weights}

A significant question in the context of our model and algorithm concerns
the structure of the attack chosen by the adversary.  The adversary is
choosing continuous values and has great leeway in how to choose them; 
potentially, for example, the adversary could choose them uniformly equal
(which, we would argue, would make the model quite uninteresting).  The 
experiments in this section address these issues.

Table \ref{histo1} describes the distribution of $x_{ij}$ values at
termination of the algorithm, for a number of networks and attack budgets.
For each test we show first (in parentheses) the number of nodes and arcs,
followed by the the attack budget and constraint set.  The data for each test
shows, for each range of resistance values, the number of arcs whose
resistance falls in that range.  

\begin{table}[ht]
\caption{{\bf Solution histogram}}
\label{histo1}
\hskip 12 pt
\centering
\begin{tabular}{|cc||cc||cc||}
\hline
\bmath{(49,90)}& \bmath{\Delta B = 57, \Gamma(3)}& \bmath{(300,409)}& \bmath{ \Delta B = 27, \Gamma(2)} &\bmath{(600,990)}& \bmath{\Delta B = 36, \Gamma(2)}\\
\hline
\hline
{\bf Range} & {\bf Count} & {\bf Range} & {\bf Count} &{\bf Range} &  {\bf Count}\\
\hline
$\left[\, 1 , \, 1 \, \right]$ & 8 &  $\left[\, 1 , \, 1  \, \right]$  & 1 & $\left[\, 1 , \, 1  \, \right]$ & 14\\
$\left(\, 1 , 2 \, \right]$ & 72 & $\left(\, 1 , \, 2  \,\right]$  & 405 & $\left(\, 1 , \, 2  \, \right]$ & 970\\
$\left(\, 2 , 3 \, \right]$ & 4 &   $\left(\, 2 , \, 9 \,  \right]$ & 0 & $\left(\, 2 , \, 5  \, \right]$ & 3\\
$\left(\, 5 , 6 \, \right]$ & 1 &  $\left(\, 9 , \, 10  \, \right]$ & 3 & $\left(\, 5 , \, 6  \, \right]$& 0\\
$\left(\, 6 , 7 \, \right]$ & 1 &    &  & $\left(\, 6 , \, 7  \, \right]$ & 1 \\
$\left(\, 7 , 8 \, \right]$ & 4 &    &  & $\left(\, 7 , \, 9  \, \right]$ & 0 \\
$\left(\, 8 , 20 \, \right]$ & 0 &   &  & $\left(\, 9 , \, 10  \, \right]$ & 2\\
\hline
\end{tabular}
\end{table}

Note that in each test case the adversary can increase the resistance of up to 
(roughly) three arcs to their maximum value. The pattern we observe
in the table is that in all three cases (i) many resistances take relatively small values
and (ii) a small number of arcs have high resistance.  Recall that 
for set $\Gamma(2)$ we always have $x^{max}_{ij} = 10$, thus in the case of
the $(300,409)$ network exactly three arcs are in the top range,
while for the $(600,990)$ network two are in the top range and one more
has relatively high resistance.  In the case of the small network there
is also a concentration 'at the top' though not in the very highest segment.
We have observed this type of behavior in many runs.  

In summary, thus, the solutions produced by the algorithm appear to superimpose
two separate effects.  From a vulnerability perspective, as we will see in the
next section, both effects play an important role.

\subsubsection{Comparison with the minimum-cardinality attack model}

The experiments in this section 
have as a first goal to effect a comparison with the $N - k$ model
as embodied by the mixed-integer programming
approach considered in Section \ref{mip}.  A direct comparison on a
case-by-case basis is not possible for a number of reasons (more on this
below) but the purpose of the tests is to investigate whether on ``similar''
data the two models behave in similar ways.

A second goal of the experiments is to investigate the impact of one of
our modeling assumptions (assumption (III) in Section \ref{nlp}), namely that demands and supplies are fixed.  Ideally,
our model should be {\em robust}, that is to say, the attack computed in a
run of the algorithm should remain effective even if the controller has the
power to adjust demands.\\

A common thread runs through both goals.  Turning to the first goal, 
it turns out that the modeling assumption (III) is, in fact, what makes a direct comparison with
the $N - k$ model difficult.  In principle, in the model in Section 
\ref{mip} one could set the desired minimum throughput to $100 \%$, 
i.e. set $T^{min} = 1.0$. But in that case an attack that disconnects a
demand node, even one with tiny demand, would be considered a success for
the attacker.  \\

To deal with these issues and still obtain a meaningful comparison, 
we set an example with $49$ nodes and $88$ arcs, and 
an example with $49$ nodes and $90$ arcs, in which 
no demand or generator node can be disconnected from the rest by removing up to three arcs.
In each case there are $4$ generators and $14$ demand nodes. A 
family of problem instances was then obtained by scaling up all capacities by a common constant.\\

In terms of the mixed-integer programming model, in each instance we
constructed a single configuration problem (generator lower bounds = 0) 
with $T^{min} = 1$, with the goal of investigating
its vulnerability should up to three arcs be removed.  Here we remind the
reader that the algorithms
\ref{mip} seek a minimum-cardinality attack that defeat the controller, and
not the most severe attack of a given cardinality.  Once our
problem is solved
the optimal attack is certified to be successful (and of minimum-cardinality), but
not necessarily {\em the} most severe attack of {\em that} cardinality.  Nevertheless,
by adjusting our formulation (\ref{mincard2})-(\ref{zero2}) we can search for
{\em a} successful attack of any given cardinality, if it exists.  The problem we obtain
is:
 \begin{eqnarray}
t^* \ \ = \ \ \mbox{ \ max} \,\, t && \label{maxcong1} \\
\mbox{Subject to:}  \ \ \ \ \ \ \ \ \ \ \ \ \ \ \ \ \  \sum_{(i,j)} z_{ij}  & \le & k, \label{foo0}\\
w_{\cC}^T \, \psi^{\cC} \, - \, t & \ge & 0, \,\,\,\,\,\, \forall \,\, \cC \subseteq \cG, \label{foo1} \\
A \psi^{\cC} \, + \, B z & \le & b \, + \, B \,\,\,\,\,\, \forall \,\, \cC \subseteq \cG, \label{foo2} \\
z_{ij} & = &  0 \,\, \mbox{or} \,\, 1, \,\,\, \forall \,\, (i,j). \label{foo3}
\end{eqnarray}
\noindent where $k$ (= 3) is a the number of arcs that the attacker can be remove.  However,
all this formulation guarantees is that $t^* > 1$ if and only if a successful attack of
cardinality $\le k$ exists -- because of the nature of our formulation, when $t^* > 1$ 
then $t^*$ will be an approximation (in general, close) to the highest severity.  
A final detail is that since $3$ lines will not disconnect the 
demands from the generators, the ``severity'' of an attack as per
formulation (\ref{foo0})-(\ref{foo3}) is
the maximum arc congestion post-attack; thus putting the problem on a common ground with the 
nonlinear models we consider. \\

\noindent For our experiments we used $\Gamma(1)$ (which allows resistances to increase by up to a factor of $20$) with an excess budget of $60$, on the network 
with 49 nodes, 90 arcs, 4 generators and 14 demand nodes.  
Note that the parameters allow the attacker to
concentrate the budget on three arcs.\\

\noindent Table \ref{comparo} contains the results.  Each row corresponds
to a different experiment, where the value indicated by \bmath{\sigma} was
used to scale all capacities (with respect to the original network).  
As $\sigma$ increases the network becomes 
progressively more difficult to interdict.

In the 'MIP' section, 
the column headed 'Cong' indicates the congestion (max. arc overload)
in the network obtained 
by removing the arcs produced by the mixed-integer programming model, 
and the column headed 'ATTACK' indicates which arcs were removed by the MIP.

In the 'NONLINEAR' section, 'Cong' indicates the maximum congestion resulting
from the increase in resistances computed by the model.  We also list the
six arcs with highest resistance (and the resistance values).  

The
column headed 'Impact' indicates the maximum congestion obtained by 
deleting the three arcs with maximum resistance (as computed by the model),
while leaving all other resistances unchanged.

We also performed additional tests with our second goal in mind, that is to
say, testing the robustness of our solutions with respect to decreased
demand levels.  In
the first test, we removed the top three (post-attack) highest resistance arcs,
while keeping all other resistances unchanged, while
allowing the controller to reduce total demand by up to $10 \%$ with the
objective of minimizing the maximum congestion.  This computation
can be formulated as a linear program; the resulting minimum congestion
value is shown in the column labeled 'I-$10\%$'.  Note that to some degree
this test also addresses the comparison with the $N - k$ model.

Similarly, but now using {\em all} resistance values as computed by the nonlinear model, and
without removing any arcs, we allowed the controller to reduce total demand by
up to $10 \%$, again with the
objective of minimizing the maximum congestion.  The column labeled 
'C-$10\%$' shows the resulting congestion value.  \\

\begin{table}[ht]
\caption{{\bf Comparison between models} } \label{comparo}
\hskip 12 pt
\centering
\begin{tabular}{|c||c|c||c|c|c|c|c|}
\hline
\bmath{\sigma}& \multicolumn{2}{|c||} {\bf MIP} & \multicolumn{5}{|c|}{{\bf NONLINEAR}} \\
\hline
  & {\bf Cong} & {\bf Attack} & {\bf Cong} & {\bf Top 6 Arcs} & {\bf Impact} &{\bf I-\bmath{10\%}} &{\bf C-\bmath{10\%}}\\
\hline
& &  & & 29(7.79), 27(7.20), 41(7.03), &  & & \\
1.0 &  1.44088 & 29,32,45 & 2.14967 & 67(7.02), 54(6.72), 79(5.71) & 1.71758 & 1.33454 & 1.67145\\
\hline
& &  & & 29(8.28), 27(7.72), 41(7.32), &  & & \\
1.2 &  1.43132 & 27,29,41 & 1.78687 &67(7.19), 54(6.92), 79(5.78)   & 1.43132 & 1.11211 & 1.38642\\
\hline
& &  & & 29(8.31), 27(7.74), 41(7.53), &  & &\\
1.4 &  1.22685 & 27,29,41 & 1.55634 & 67(7.48), 54(7.18), 79(6.15)  & 1.22685 & 0.95324 & 1.21329\\
\hline
& &  & & 29(8.18), 27(7.58), 41(7.53), &  & &\\
1.6 &  1.07349 & 27,29,41 & 1.35995 & 67(7.58), 54(7.22), 79(6.25)  & 1.07349 & 0.83409 & 1.05458\\
\hline
& &  & & 29(8.43), 27(7.90), 41(7.53), &  & &\\
1.8 &  0.692489 & 18,57,60 & 1.20271 &  67(7.48), 54(7.18), 79(6.12)  & 0.95421 & 0.74141 & 0.93595\\
\hline
& &  & &  29(7.87), 27(7.29), 41(7.04),  & & &\\
2.0 & 0.68630 & 20,89,45 & 1.07733 & 67(7.01), 54(6.70), 79(5.63)&  0.85889  & 0.66727 & 0.83878\\
\hline
\end{tabular}
\end{table}

\noindent {\bf Comments.} As before,
we see that the
solutions to the nonlinear model tend to concentrate the attack on a relatively small number of lines, while
at the same time investing small portions of the attack budget on other lines.
This helps highlight the significant 
overlap between
the results from the two models.  Note that in the cases for 
$\sigma \ = \ 1.2, \, 1.4, \, 1.6$ the set of attacked lines show high
correlation.

Moreover, the two models are consistent: the severity of the attack
as measured by the maximum congestion levels (the 'Cong' parameters),
for both models,  
decrease as the scale increases (as one should expect). 

The last three columns of the table address our second set of questions --
they appear to show that the solution computed by the nonlinear model is
robust; even as the controller reduces total demand, the congestion level 
is proportionally reduced.  Finally, note that the congestion values in the 'Impact' column 
are significantly smaller than the corresponding values in the 'Cong'
column; similarly, the 'C-$10\%$' values are higher than the
'I-$10\%$' values -- thus, the low $x_{ij}$ arcs in the nonlinear attack {\em do}
play a significant role.

\subsection{Future work}
We find our experiments with the nonlinear model highly encouraging.  A 
follow-up project that we are now starting is to embed the above approach
in a global optimization procedure.  A related project that we are also
undertaking in parallel is the extension of the above ideas to more
complex models of power flows.

\end{document}